\documentclass[doublespacing]{elsart}
\usepackage{graphicx}
\usepackage{subfigure}
\usepackage{epsfig}
\usepackage{amsmath}
\usepackage{amsfonts}   
\usepackage{amssymb}
\journal{}

\newcommand{\mlabel}[1]{\label{#1}
}

\newcommand{\seq}{\begin{equation}}                 
\newcommand{\eeq}[1]{\label{#1}\end{equation}
    }

\newcommand{\epf}{$ \quad \Box$ \par \vspace{1ex}}

\newtheorem{Theorem}{Theorem}[section]
\newcommand{\sthm}{\begin{Theorem}}         
\newcommand{\ethm}{\end{Theorem}}           

\newtheorem{Corollary}[Theorem]{Corollary}
\newcommand{\scor}{\begin{Corollary}}       
\newcommand{\ecor}{\end{Corollary}}         

\newtheorem{Lemma}[Theorem]{Lemma}
\newcommand{\slm}{\begin{Lemma}}            
\newcommand{\elm}{\end{Lemma}}              

\newtheorem{Proposition}[Theorem]{Proposition}
\newcommand{\spro}{\begin{Proposition}}            
\newcommand{\epro}{\end{Proposition}}              

\newtheorem{Example}[Theorem]{\sc Example}
\newcommand{\sex}{\begin{Example}\rm}        
\newcommand{\eex}{\end{Example}}             

\newtheorem{Remark}[Theorem]{Remark}
\newcommand{\srmark}{\begin{Remark}\rm}        
\newcommand{\ermark}{\end{Remark}}             

\newtheorem{Definition}[Theorem]{Definition}
\newcommand{\strdef}{\begin{Definition}\rm}        
\newcommand{\eeddef}{\end{Definition}}             

\newcommand{\seql}{\begin{eqnarray*}}       
\newcommand{\eeql}{\end{eqnarray*}}

\newcommand{\smlist}[1]{\begin{list}           
                      {(#1{zzcount})}{\usecounter{zzcount}}}
\newcommand{\elist}{\end{list}}

\newcommand{\vect}[1]{\mathbf{#1}}
\providecommand{\abs}[1]{\lvert#1\rvert}
\providecommand{\norm}[1]{\lVert#1\rVert}

\begin{document}
\begin{frontmatter}

\title{Spreading speeds and traveling waves for non-cooperative reaction-diffusion systems\thanksref{ref1}}
\author{Haiyan Wang}
\address{Division of Mathematical and Natural Sciences\\ Arizona State University \\ Phoenix, AZ 85069-7100, USA}
\thanks[ref1]{This manuscript is available at arXiv:1007.0950}
\ead{wangh@asu.edu}
\begin{abstract}
Much has been studied on the spreading speed and traveling wave solutions for cooperative reaction-diffusion systems. In this paper, we shall
establish the spreading speed for a large class of non-cooperative reaction-diffusion systems and characterize the spreading speed as the slowest
speed of a family of non-constant traveling wave solutions. As an application, our results are applied to a partially cooperative system describing
interactions between ungulates and grass.
\end{abstract}

\end{frontmatter}

\pagenumbering{arabic}

\section{Introduction}\label{induc}

Fisher \cite{Fisher} studied the nonlinear parabolic equation
\begin{equation}\label{eq0}
w_t=w_{xx}+w(1-w).
\end{equation}
for the spatial spread of an advantageous gene in a population and conjectured $c^*$ is the asymptotic speed of propagation of the advantageous gene.
His results show that (\ref{eq0}) has a traveling wave solution of the form $w(x+ct)$ if only if $|c| \geq c^*=2.$ Kolmogorov, Petrowski, and Piscounov
\cite{Kolmogorov} proved the similar results with more general model. Those pioneering work along with the paper by
Aronson and Weinberger \cite{Aronson1975,Aronson1978} confirmed the conjecture of Fisher and established the speeding spreads for nonlinear parabolic equations.
Lui \cite{Lui1989} established the theory of spreading speeds for cooperative recursion systems. In a series of papers,
Weinberger, Lewis and  Li \cite{Weinberger2002,Weinberger2005,Weinberger2002-1,Weinberger2007} studied spreading speeds and traveling waves for more general
cooperative recursion systems, and in particular, for quite general cooperative reaction-diffusion systems by analyzing of
traveling waves and the convergence of initial data to wave solutions.

However, due to various biological or physical constrains, many reaction-diffusion systems are not necessarily cooperative.
Thieme \cite{Thieme1979} showed that asymptotic spreading speed of a model with
nonmonotone growth functions can still be obtained by constructing monotone functions. Weinberger, Kawasaki and Shigesada \cite{WeinbergerKS2009} discussed
the reaction-diffusion model
\begin{equation}\label{eqoper0}
\begin{split}
\frac{\partial u_1}{\partial t}&=d_1 \Delta u_1+u_1[-\alpha -\delta u_1+r_1 u_2]\\
\frac{\partial u_2}{\partial t}&=d_2 \Delta u_2+u_2 r_2[1-u_2 +h(u_1)]\\
\end{split}
\end{equation}
where $ d_1, \alpha, \delta, r_1, d_2, r_2$ are all positive parameters. This system describes the interaction between ungulates with
linear density $u_1(x, t)$ and grass with linear density $u_2(x, t).$  The function $h(u_1)$ models the increase in the specific growth
rate of the grass due to the presence of ungulates. When the density $u_1$ is small the net effect of ungulates
is increasingly beneficial, but as the density increases above a certain value, the benefits decrease with increasing $u_1$.
(\ref{eqoper0}) is partially cooperative 2-species reaction-diffusion model, meaning that it is cooperative for small population densities
but not for large ones. By employing comparison methods \cite{WeinbergerKS2009} established spreading speeds for (\ref{eqoper0}).
In Section \ref{example}, we take the non-monotone Ricker function $u_1e^{-u_1}$ as $h(u_1)$, which is simpler than that of \cite{WeinbergerKS2009},
and apply our main theorem (Theorem \ref{th30}) to (\ref{eqoper0}). The application of our general theorem allows us to
characterize the spreading speed as the slowest speed of traveling wave solutions to (\ref{eqoper0}), which is new and was not
proved in \cite{WeinbergerKS2009}.  Non-cooperative reaction-diffusion systems frequently occur in other biological systems such as epidermal wound healing
(see, Sherratt and Murray \cite{Sherratt1990,Sherratt1991}, Dale, Maini, Sherratt \cite{Dale1994}). In a recent paper by the author \cite{HwangSystemPDE2}, spreading speeds and traveling waves
for a non-cooperative reaction-diffusion model of epidermal wound healing were established.

For related non-monotone integro-difference equations, Hsu and Zhao \cite{hsu2008}, Li, Lewis and Weinberger \cite{LiLewis2009}
extended the theory of spreading speed and established the existence of travel wave solutions.
The author and Castillo-Chavez \cite{HwangIntegralDiff} prove that a large class of nonmonotone integro-difference systems have spreading speeds and traveling wave solutions.
Such an extension is largely based on the construction of two monotone operators with appropriate properties and fixed point theorems in Banach spaces.
A similar method was also used in Ma \cite{ma2007} and the author \cite{Hwang2009} to
prove the existence of traveling wave solutions of nonmonotone reaction-diffusion equations.

In this paper, we shall establish the spreading speed for a general non-cooperative system (\ref{eq1}) and characterize its spread speed
as the slowest speed of a family of non-constant traveling wave solutions of (\ref{eq1}). Our main theorem (Theorem \ref{th30}) will be applied
to (\ref{eqoper0}) in Section \ref{example}.

We begin with some notation. We shall use $R, k, k^{\pm}, f, f^{\pm}, r, u,v$ to denote vectors in $\mathbb{R}^N$ or $N$-vector valued functions
, and $x,y,\xi$ the single variable in $\mathbb{R}$. Let $u=(u_i), v=(v_i) \in \mathbb{R}^N$, we write $u \geq v$ if $u_i \geq v_i $ for all $i$;
and
$u \gg v$ if $u_i > v_i$ for all $i$. A vector $u$ is positive if $u \gg 0$. For any $r=(r_i) \gg 0, r \in \mathbb{R}^N$ let
$$
[0,r]= \{ u: 0 \leq u \leq r, u \in \mathbb{R}^N\}\subseteq \mathbb{R}^N
$$
and
$$
\mathcal{C}_{r}= \{u=(u_i): u_i \in C(\mathbb{R}, \mathbb{R}), 0\leq u_i(x) \leq r_i,  x\in \mathbb{R},\; i=1,...,N\},
$$
where $C(\mathbb{R}, \mathbb{R})$ is the set of all continuous functions from $\mathbb{R}$ to $\mathbb{R}$.

Consider the system of reaction-diffusion equations
\begin{equation}\label{eq1}
u_t=Du_{xx}+f(u), x \in \mathbb{R},\; t\geq 0.
\end{equation}
with
\begin{equation}\label{eq1bc}
u(x,0)=u_0(x),\;\;  x \in \mathbb{R},
\end{equation}
where $u=(u_i)$, $D=\text{diag} (d_1, d_2, ...,d_N), d_i>0, i=1,...,N$
$$f(u)=(f_1(u),f_2(u),...,f_N(u)),$$
$u_0(x)$ is a bounded uniformly continuous function on $\mathbb{R}.$ In this paper, by a solution we mean a twice continuously differentiable function $u$
satisfying appropriate equation in $\mathbb{R}$ and an initial condition.

In order to deal with non-cooperative system, we shall assume that there are additional
two monotone operators $f^{\pm}$, one lies above and another below $f$ with the corresponding equations
 \begin{equation}\label{eq1+}
u_t=Du_{xx}+f^+(u), x \in \mathbb{R},\; t\geq 0.
\end{equation}
\begin{equation}\label{eq1-}
u_t=Du_{xx}+f^-(u), x \in \mathbb{R},\; t\geq 0.
\end{equation}
Such an assumption will enable us to make use of the corresponding results for cooperative systems in
\cite{Lui1989, Weinberger2002-1} to establish spreading speeds for (\ref{eq1}).

\begin{enumerate}
    \item[(H1)] \begin{itemize}
             \item[i.] Let $ k^{+}=(k^+_i)>>0$ and $f: [0,k^{+}] \to \mathbb{R}^N$ be a continuous and twice piecewise continuous differentiable function.
             Assume that $\mathcal{C}_{k^{+}}$ is an invariant set of (\ref{eq1}) in the sense that
             for any given $u_0 \in \mathcal{C}_{k^+}$, the solution of (\ref{eq1}) with the initial condition $u_0$ exists and remains
             in $\mathcal{C}_{k^{+}}$ for $t \in [0,\infty)$.
             \item[ii.] Let $ 0<<k^{-}=(k^-_i)\leq k=(k_i) \leq k^{+}.$
              Assume there exist continuous and twice piecewise continuous differentiable function $f^{\pm}=(f^{\pm}_i): [0,k^+] \to \mathbb{R}^N$
              such that for $u \in [0,k^+]$
              \begin{equation*}
              f^{-}(u) \leq f(u) \leq f^{+}(u).
              \end{equation*}
              \item[iii.] $f(0)=f(k)=0$ and there is no other positive equilibrium of $f$ between $0$ and $k$.
              $f^{\pm}(0)=f^{\pm}(k^{\pm})=0$. There is no other positive equilibrium of $f^{\pm}$ between $0$ and $k^{\pm}$.
              \item[iv.] (\ref{eq1+}) and (\ref{eq1-}) are cooperative
              (i.e. $\partial_if^{\pm}_j(u) \geq 0$ for $u \in [0, k^{\pm}], i\neq j$). $f^{\pm}(u), f(u)$ have the same Jacobian matrix $f'(0)$ at $u=0$.
              \end{itemize}
\end{enumerate}

A traveling wave solution $u$ of (\ref{eq1}) is  a solution of the form $u=u(x+ct), u \in C(\mathbb{R}, \mathbb{R}^N)$ .
Substituting $u(x,t)=u(x+ct)$ into (\ref{eq1}) and letting $\xi=x+ct$, we obtain the wave equation
\begin{equation}\label{eq211}
Du''(\xi)-cu'(\xi)+f(u(\xi))=0, \;\; \xi \in \mathbb{R}.
\end{equation}
Now if we look for a solution of the form
$(u_i)=\big ( e^{\lambda \xi}\eta^i_{\lambda}\big),\lambda>0, \eta_{\lambda}=(\eta^i_{\lambda})>>0$ for the linearization of (\ref{eq211}) at the origin,
we arrive at the following system equation
\begin{equation*}
\text{diag}(d_i \lambda^2 -c \lambda)\eta_{\lambda}+f'(0)\eta_{\lambda}=0
\end{equation*}
which can be rewritten as the following eigenvalue problem
\begin{equation}\label{egenvalue}
\frac{1}{\lambda}A_{\lambda}\vect{\eta_{\lambda}}= c \vect{\eta_{\lambda}},
\end{equation}
where
\begin{equation*}
A_{\lambda}=(a^{i,j}_{\lambda})=\text{diag}(d_i\lambda^2)+f'(0)
\end{equation*}

The matrix $f'(0)$ has nonnegative off diagonal elements. In fact, there is a constant $\alpha$ such that $f'(0)+\alpha I$ has
nonnegative entries, where $I$ is the identity matrix.

By reordering the coordinates, we can assume that $f'(0)$ is in block lower triangular
form, in which all the diagonal blocks are irreducible or $1$ by $1$ zero matrix. A matrix is irreducible if it is not similar
to a lower triangular block matrix with two blocks via a permutation. From the Perron-Frobenius theorem any irreducible
matrix $A$ with nonnegative entries has a
unique principal positive eigenvalue( which is the spectral radius  of $A$, $\rho(A)$)  with a corresponding principal eigenvector
with strictly positive coordinates.  For an irreducible matrix $A$ with nonnegative off diagonal elements,
we shall call the eigenvalue $\rho(A+\alpha I)-\alpha$ of $A$, which has the same
positive eigenvector, the principal eigenvalue of $A$ (see e.g. \cite{Horn,Weinberger2002-1}).
Let $$\Psi(A)=\rho(A+\alpha I)-\alpha.$$
Here $A+\alpha I$ is irreducible and nonnegative, and $\rho(A+\alpha I)$ is the spectral radius of $A+\alpha I.$

We shall need the following assumption (H2). Notice that (H2) is assumed for $\lambda=0$ in \cite{Weinberger2002-1}.
However, with (H2), we are able to obtain better estimates for traveling solutions and the minimum speed $c^*$, see Lemma \ref{lmeigen}.
As a result, for the example in Section \ref{example}, a slightly stronger assumption ($d_1 \geq d_2$) than \cite{WeinbergerKS2009} is be
assumed.

\begin{enumerate}
  \item[(H2)] \begin{itemize}
   \item[] Assume that, for each $\lambda>0$, $A_{\lambda}$ is in block lower triangular
form and the first diagonal block has the positive principal eigenvalue  $\Psi(A_{\lambda})$, and $\Psi(A_{\lambda})$ is strictly larger than the principal
eigenvalues of all other diagonal blocks, and that there is a positive eigenvector $\nu_{\lambda}=(\nu^i_{\lambda})>>0$
of  $A_{\lambda}$ corresponding to  $\Psi(A_{\lambda})$. And further assume that $\nu_{\lambda}$ is continuous with respect to $\lambda.$
   \end{itemize}
\end{enumerate}

Let $$
\Phi(\lambda)=\frac{1}{\lambda} \Psi(A_{\lambda})> 0.
$$
\begin{figure}
\begin{center}
  \includegraphics[width=7cm]{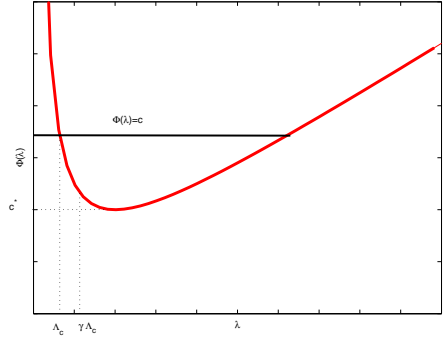}
  \caption{The red curve is $\Phi(\lambda)$. The minimum of $\Phi(\lambda)$ is $c^*$. For $c>c^*$, the left solution of $\Phi(\lambda)=c$ is $\Lambda_c$.} \label{fig4}
\end{center}
\end{figure}
According to Lemma \ref{lmeigen}, we can expect the graph of $\Phi$ as in Fig. \ref{fig4}. For the example in Section \ref{example},
$\Phi$ is a strictly convex function of $\lambda$ and, clearly satisfies Lemma \ref{lmeigen}.

Now we state Lemma \ref{lmeigen}, which shall enable us to give accurate asymptotic estimates of traveling solutions.
Lemma \ref{lmeigen} is a analogous result in Weinberger \cite{Weinberger1978} and Lui \cite{Lui1989}.
However, due to the fact that $f'(0)$ is only quasi-positive and the elements of $A_{\lambda}$ are not necessarily log convex,
some of its proofs here are different from Lui \cite{Lui1989}.  A theorem on the convexity of the dominant eigenvalue of matrices
due to Cohen \cite{Cohen1981} is used to show that $\Psi(A_{\lambda})$ is convex function of $\lambda$.

There are two direct consequences  of Lemma \ref{lmeigen}. First one, it improves \cite[Theorem 4.2]{Weinberger2002-1}, which will be used in this paper,
by eliminating the case (b) in \cite[Theorem 4.2]{Weinberger2002-1} because of Lemma \ref{lmeigen} (7). Second, Lemma \ref{lmeigen} (8)
will allow us to construct explicit lower solutions and therefore asymptotic behavior of traveling solutions of (\ref{eq1}) can be obtained.

\slm\mlabel{lmeigen}  Assume that $(H1)-(H2)$ hold. Then
\begin{enumerate}
  \item [(1)] $\Phi(\lambda) \to \infty$ as $\lambda \to 0;$
  \item [(2)]  $\Phi(\lambda) \to \infty$ as $\lambda \to \infty;$
  \item [(3)] $\Phi(\lambda)$ is decreasing as $\lambda=0^+;$
  \item [(4)]  $\Psi(A_{\lambda})$ is a convex function of $\lambda>0;$
  \item [(6)] $\Phi'(\lambda) $ changes sign at most once on $(0, \infty)$
  \item [(7)] $\Phi(\lambda)$ has the minimum $$
c^*=\inf_{\lambda>0}\Phi(\lambda)>0
$$
at a finite $\lambda$.
  \item [(8)] For each $c > c^*$, there exist $\Lambda_{c}>0$ and $\gamma \in (1,2)$ such that
$$\Phi(\Lambda_{c})=c,\;\; \Phi(\gamma \Lambda_{c})<c.$$
That is $$\frac{1}{\Lambda_c}A_{\Lambda_c} \nu_{\Lambda_c}= \Phi(\Lambda_{c})\nu_{\Lambda_c}=c\nu_{\Lambda_c}$$
and
$$\frac{1}{\gamma\Lambda_c}A_{\gamma\Lambda_c} \nu_{\gamma\Lambda_c}= \Phi(\gamma\Lambda_{c})\nu_{\gamma\Lambda_c}<c\nu_{\gamma\Lambda_c}$$
where $\nu_{\Lambda_c}>>0,\nu_{\gamma\Lambda_c}>>0$ are positive eigenvectors
of $\frac{1}{\Lambda_c}A_{\Lambda_c}, \frac{1}{\gamma\Lambda_c}A_{\gamma\Lambda_c}$
corresponding to eigenvalues $\Phi(\Lambda_{c})$ and $\Phi(\gamma\Lambda_{c})$ respectively.
\end{enumerate}
\elm
\pf
We only need to prove those different from \cite{Weinberger1978,Lui1989}. The proof of the convexity of $\Psi(A_{\lambda})$ is
similar to that in  Crooks \cite{Crooks1996} for matrices with positive off-diagonal elements. It is easily seen that
$\Psi(A_{\lambda})=\rho(A_{\lambda}+\alpha I)-\alpha$ is non-decreasing function of $\lambda>0$ (\cite[Theorem 8.1.18]{Horn}).
Further, a theorem on the convexity of the dominant eigenvalue of matrices due to Cohen \cite{Cohen1981} states that for
any positive diagonal matrices $D_1, D_2$ and $t \in (0,1)$,
$$
\Psi(tD_1+(1+t)D_2+f'(0)) \leq t\Psi(D_1+f'(0))+(1-t)\Psi(D_2+f'(0))
$$
As before $\Psi(A)$ is the principle eigenvalue of $A$. Now if $\alpha_1, \alpha_2 \in \mathbb{R}$ and $t \in (0,1)$,
$$(t\alpha_1+(1-t)\alpha_2)^2 \leq t\alpha_1^2+(1-t)\alpha_2^2.$$
This implies that
\begin{equation*}
\begin{split}
\Psi(A_{t\lambda_1+(1-t)\lambda_2})& = \Psi((t\lambda_1+(1-t)\lambda_2)^2D+f'(0))\\
                       &\leq \Psi(t\lambda_1^2D+(1-t)\lambda_2^2D+f'(0))\\
                       &\leq t\Psi(\lambda_1^2D+f'(0))\\
                       & \quad +(1-t)\Psi(\lambda_2^2D+f'(0))\\
                       &=t\Psi(A_{\lambda_1})+(1-t)\Psi(A_{\lambda_2})\\
\end{split}
\end{equation*}
Since $\Psi(A_{\lambda})$ is a simple root of the characteristic equation of an irreducible block,
it can be shown that $\Psi(A_{\lambda})$ is twice continuously differentiable on $\mathbb{R}$. Thus
$$\Psi''(\lambda)\geq 0$$
and a calculation shows
$$
[\lambda\Phi(\lambda)]'=\Psi'(\lambda)
$$
$$
\Phi'(\lambda)=\frac{1}{\lambda}[\Psi'(\lambda)-\Phi(\lambda)]
$$
and
$$
(\lambda^2\Phi'(\lambda))'=\lambda\Psi''(\lambda)\geq 0.
$$
(6) is a consequence of the above inequalities.  As for (2), we  need to prove that
$\lim_{\lambda \to \infty}\frac{\Psi(A_{\lambda})}{\lambda}=\infty$. In fact, there exists an $\epsilon>0$ such that all diagonal elements of
$D-\epsilon I$ are strictly positive, then $\Psi\big(D-\epsilon I\big)>0$ and choose $\lambda$ large enough so that
\begin{equation*}
\begin{split}
\Psi(A_{\lambda})& =\Psi(D\lambda^2+f'(0)) \\
                       &=\Psi\big((D-\epsilon I)\lambda^2+ (\epsilon \lambda^2 I+f'(0))\big)\\
                       &\geq \Psi\big((D-\epsilon I)\lambda^2\big)\\
                       & = \lambda^2 \Psi\big(D-\epsilon I\big)
\end{split}
\end{equation*}
Thus $\lim_{\lambda \to \infty}\frac{\Psi(A_{\lambda})}{\lambda}=\infty$.  As we discussed before, (H2) implies the existence of positive
eigenvector $\nu_{\lambda}>>0$ corresponding to $\Phi(A_{\lambda})$.
The first statement of (8) is a consequence of (1)-(7). The second statement of (8) is just a rephrase
of the fact that $\nu_{\lambda}>>0$ is a eigenvector of $\frac{1}{\lambda}A_{\lambda}$ corresponding to eigenvalue $\Phi(A_{\lambda})$
for $\lambda=\Lambda_c$ and $\gamma \Lambda_c$.

\epf

In addition to (H1-H2), we also need assumption (H3) which only requires the nonlinearity grow less than its linearization along the particular function $\nu_{\lambda} e^{-\lambda x }$ \cite{Weinberger2002-1}. Such a condition can be satisfied for many biological systems.
\begin{enumerate}
  \item[(H3)] \begin{itemize}
   \item[] Assume that for any $\alpha>0, \lambda>0$ $$
   f^{\pm}(\alpha \nu_{\lambda}) \leq \alpha f'(0)\nu_{\lambda},\;\; \text{where} \;\; \nu_{\lambda}= (\nu^i_{\lambda}).
   $$
   \end{itemize}
\end{enumerate}

We now recall results on the spreading speeds in  Weinberger, Lewis and Li \cite{Weinberger2002-1} and Lui \cite{Lui1989}.
While Theorem 4.1 \cite{Weinberger2002-1} holds for non cooperative reaction-diffusion systems, it does require that the reaction-diffusion system has a single speed.
In general, such a condition is very difficult to verify. In the same section, for cooperative systems, Theorem 4.2 in \cite{Weinberger2002-1}
provides sufficient conditions to have a single speed. The following theorem combines the results of Theorems 4.1 and 4.2 in \cite{Weinberger2002-1},
which can be a consequence of Theorems 3.1 and 3.2 for discrete-time recursions in Lui \cite{Lui1989}.

\sthm (Weinberger, Lewis and Li \cite{Weinberger2002-1})\label{th20} Assume $(H1)-(H3)$ hold and (\ref{eq1}) is \textbf{cooperative}.  Then the following statements are valid:
\begin{enumerate}
\item [(i)] For any $u_0 \in \mathcal{C}_{k}$ with compact support, let $u(x,t)$ be the solution of (\ref{eq1}) with (\ref{eq1bc}).
Then
$$
\lim_{t \to \infty}  \sup_{\abs{x}\geq ct} u(x,t) = 0, \text{ for } c> c^*
$$

\item [(ii)] For any strictly positive vector $\omega \in \mathbb{R}^N$, there is a positive $R_{\omega}$ with the property that if $u_0 \in \mathcal{C}_{k}$
and $u_0 \geq \omega $
on an interval of length $2R_{\omega}$, then the solution $u(x,t)$ of (\ref{eq1})  with (\ref{eq1bc}) satisfies
$$
\liminf_{t \to \infty} \inf_{\abs{x}\leq tc} u(x,t) =  k, \text{ for } 0< c< c^*
$$
\end{enumerate}
\ethm
In another paper \cite{Weinberger2005}, for \textbf{cooperative systems}, Li, Weinberger and Lewis established that the slowest spreading
speed $c^*$ can always be characterized as the slowest speed of a family of traveling waves.
These results describe the properties of spreading speed $c^*$ for monotone systems. Based on these spreading results for cooperative systems,
we will discuss analogous spreading speed results for non cooperative systems.

\section{Main Results}\label{rest}
Our new contributions in this paper are to establish the spreading speed
(Theorem \ref{th30} (i-ii)) for general non-cooperative reaction-diffusion systems (\ref{eq1}),
and further characterize the spreading speed as the speed of the slowest non-constant
traveling wave solutions (Theorem \ref{th30} (iii-v)).

Although the existence of traveling wave solutions for cooperative systems are known (see,e.g. \cite{Weinberger2005}), we shall prove the existence of
traveling wave solutions for both cooperative and non cooperative systems as our proofs for non cooperative systems are based on those for cooperative systems.
Further, in additions to the existence of traveling wave solutions, we shall be able to obtain asymptotic behavior of
the traveling wave solutions in terms of eigenvalues and eigenvectors for both cooperative and non cooperative systems.
The following theorem is our main results.

\sthm\label{th30} Assume $(H1)-(H3)$ hold.  Then the following statements are valid:
\begin{enumerate}
\item [(i.)] For any $u_0 \in \mathcal{C}_{k}$ with compact support, the solution $u(x,t)$ of (\ref{eq1})  with (\ref{eq1bc}) satisfies
$$
\lim_{t \to \infty}  \sup_{\abs{x}\geq tc} u(x,t) = 0, \text{ for } c> c^*
$$

\item [(ii.)] For any vector $\omega \in \mathbb{R}^N, \omega>>0$, there is a positive $R_{\omega}$ with the property that if $u_0 \in \mathcal{C}_{k}$
and $u_0 \geq \omega $
on an interval of length $2R_{\omega}$, then the solution $u(x,t)$ of (\ref{eq1})  with (\ref{eq1bc}) satisfies
$$
k^- \leq \liminf_{t \to \infty} \inf_{\abs{x}\leq tc} u(x,t) \leq  k^+, \text{ for } 0< c< c^*
$$

\item [(iii.)] For each $c > c^*$ (\ref{eq1}) admits a traveling wave solution $u=u(x+ct)$ such that
$0 << u(\xi) \leq k^+, \xi \in \mathbb{R}$,
$$k^-\leq \liminf_{\xi \to \infty}u(\xi) \leq \limsup_{\xi \to \infty}u(\xi)\leq k^+$$
and
\begin{equation}\label{aystomth112}
\lim_{\xi \to -\infty} u(\xi)e^{-\Lambda_{c} \xi}=\nu_{\Lambda_{c}}.
\end{equation}
If, in addition, (\ref{eq1}) is cooperative in $\mathcal{C}_{k}$, then $u$ is nondecreasing on $\mathbb{R}$.

\item [(iv.)] For $c = c^*$ (\ref{eq1}) admits a nonconstant traveling wave solution $u=u(x+ct)$ such that
$0 \leq  u(\xi) \leq k^+, \xi \in \mathbb{R}$,
$$k^-\leq \liminf_{\xi \to \infty}u(\xi) \leq \limsup_{\xi \to \infty}u(\xi)\leq k^+.$$

\item [(v.)] For $0<c<c^*$ (\ref{eq1}) does not admit a traveling wave solution $u=u(x+ct)$
 with $ \liminf_{\xi \to \infty} u(\xi)>>0$ and $u(-\infty)=0.$
\end{enumerate}
\ethm

\srmark\label{rem1}
In many cases, $f^{\pm}$ can be taken as piecewise functions consisting of $f$ and appropriate constants as
demonstrated in Section \ref{example}. In order to have a better estimate
for the traveling wave solution $u$ for non cooperative systems, it is desirable to choose two function $f^{\pm}$ which are close enough.
The smallest monotone function above $f$ and the largest monotone function below $f$  are natural choices of $f^{\pm}$
if they satisfy other requirements, See \cite{Thieme1979,hsu2008,LiLewis2009} for the discussion for scalar cases and \cite{WeinbergerKS2009}
for a partially cooperative reaction-diffusion system. Our construct of $f^-$ in Section \ref{example} is different from the previous papers.
\ermark

\srmark\label{rem3}
The invariant set of (H1) (i) can often be established by Comparison Principle \ref{comparison1}.  In fact,
for a given $u_0 \in \mathcal{C}_{k^+}$,
let $u(x,t)$ be the solution of (\ref{eq1})
with the initial condition $u_0$.
If we can choose appropriate
$f^-, f^+$  so that $f^-(u) \leq f(u) \leq f^+(u)$ for all $u \in \mathbb{R}^N$,
it follows that
$$
k^+_t-Dk^+_{xx} -f^+(k^+)=0=u_t-Du_{xx} -f(u) \geq  u_t-Du_{xx} -f^+(u).
$$
and
$$
0-D0 -f^-(0)=0=u_t-Du_{xx} -f(u) \leq  u_t-Du_{xx} -f^-(u).
$$
Comparison Principle \ref{comparison1} implies that
$$
0 \leq u(x,t) \leq k^+, x \in \mathbb{R}, t >0.
$$
Now according to Smoller \cite[Theorem 14.4]{Smoller1994} (\ref{eq1}) (and also (\ref{eq1+}), (\ref{eq1-})) has a solution $u$ for $t \in [0, \infty)$ and $0 \leq u \leq k^+$ if
the initial value $u_0$ is uniformly continuous on $\mathbb{R}$. Now we can establish an invariance set for (\ref{eqoper00}) in Section \ref{example}.
First we can extend $h, h^+$ in (\ref{eqoper00}) and (\ref{eqoper1}) to zero for $w_1<0$. Now
let $f^+$ be the reaction terms in (\ref{eqoper1}).   And further let $f^-$ be the reaction terms in (\ref{eqoper00})
with $h$ being replaced by the constant zero function for all $w_1 \in \mathbb{R}$. From the above discussion, we can see that $[0,k^+]$ is
an invariance set for (\ref{eqoper00}).
\ermark

\srmark\label{rem5}
When (\ref{eq1}) is cooperative in $\mathcal{C}_{k}$, then $f^{\pm}=f.$
\ermark

\srmark\label{rem6}
As indicated in \cite{Weinberger2002-1}, if $f$ is not defined everywhere, (H3) can be replaced by (H3').
Without extra assumptions, (H3') can be verified with a slightly complicated computation.
\begin{enumerate}
  \item[H3'] \begin{itemize}
   \item[] For each $\lambda>0$, let $v^{\pm}=(\min\{k^{\pm}_i, \alpha \nu^i_{\lambda}\}).$ Assume that for any $\alpha>0$ $$
   f^{\pm}(v^{\pm}(x)) \leq \alpha f'(0)(\nu^i_{\lambda}).
   $$
   \end{itemize}
\end{enumerate}
\ermark

The paper is organized as follows. Theorem \ref{th30} (i)-(ii) shall be proved in Section \ref{spreadings} and  Theorem \ref{th30} (iii)-(v) in 
Section \ref{travel}.

\section{The Spreading Speed}\label{spreadings}

\subsection{Comparison Principle}\label{comparsion}
We state the following comparison theorem for cooperative systems of reaction-diffusion equations in  Weinberger, Kawasaki and Shigesada \cite{WeinbergerKS2009} .
The comparison principle is a consequence of the maximum principle (see, e.g., Protter and Weingberger \cite{Protter1984}).

%

\sthm \label{comparison1} Let $D$ be a positive definite diagonal matrix. Assume $F=(F_j)$ is vector-valued functions in $\mathbb{R}^N$ are continuous and piecewise continuously differentiable
in $\mathbb{R}$ and the underling system is cooperative in the sense that for each $j$, $F_j$ is nondecreasing
in all but the $j$th component. Suppose that  $u(x,t), v(x,t)$ satisfy
\begin{equation}\label{comp1}
\begin{split}
u_t-Du_{xx} -F(u) \leq  v_t-Dv_{xx} -F(v)
\end{split}
\end{equation}
If $u(x,t_0) \leq v(x,t_0), \;\; x\in \mathbb{R}$, then
$$
u(x,t) \leq v(x,t), \;\; x\in \mathbb{R}, t \geq t_0.
$$
\ethm


We are now able to prove Parts (i) and (ii) of Theorem \ref{th30}.

\subsection{Proof of Parts (i) and (ii) of Theorem \ref{th30}}

Part (i).  For a given $u_0 \in \mathcal{C}_{k}$ with compact support, let $u^+(x,t)$ be the solutions of (\ref{eq1+}) with the same initial
condition $u_0$, then Comparison Principle \ref{comparison1} implies that $u^+(x,t) \in \mathcal{C}_{k^+}$ and

$$
u^+_t-Du^+_{xx} -f^+(u^+)=0=u_t-Du_{xx} -f(u) \geq  u_t-Du_{xx} -f^+(u).
$$
and
$$
0-D0 -f^-(0)=0=u_t-Du_{xx} -f(u) \leq  u_t-Du_{xx} -f^-(u).
$$
Comparison Principle \ref{comparison1} further implies that
$$
0 \leq u(x,t) \leq u^{+}(x,t), x \in \mathbb{R}, t >0.
$$
Thus for any $c>c^*$, it follows from Theorem \ref{th20} (i) that $$
\lim_{t \to \infty} \sup_{\abs{x}\geq tc} u^{+}(x,t) = 0,
$$
and hence
$$
\lim_{t \to \infty} \sup_{\abs{x}\geq tc} u(x,t) = 0,
$$
Part (ii). According to  Theorem \ref{th20} (ii), for any strictly positive constant $\omega$, there is a positive $R_{\omega}$
(choose the larger one between the $R_{\omega}$ for (\ref{eq1+}) and the $R_{\omega}$ for (\ref{eq1-})) with the property that if $u_0 \geq \omega $
on an interval of length $2R_{\omega}$, then the solutions $u^{\pm}(x,t)$ of (\ref{eq1+}) and (\ref{eq1-}) with the same initial value $u_0$
are in $\mathcal{C}_{k^+}$ and satisfy
$$
\liminf_{t \to \infty} \inf_{\abs{x}\leq tc} u^{\pm}(x,t) =  k^{\pm}, \text{ for } 0< c< c^*.
$$
As before we have

$$
u^+_t-Du^+_{xx} -f^+(u^+)=0=u_t-Du_{xx} -f(u) \geq  u_t-Du_{xx} -f^+(u)
$$
and
$$
u^--Du^-_{xx} -f^-(u^-)=0=u_t-Du_{xx} -f(u) \leq  u_t-Du_{xx} -f^-(u).
$$
Thus, Comparison Principle \ref{comparison1} implies that

$$
u^-(x,t) \leq u(x,t) \leq u^{+}(x,t), x \in \mathbb{R}, t >0.
$$
Thus for any $c<c^*$, it follow from Theorem \ref{th20} (ii) that $$
\liminf_{t \to \infty} \inf_{\abs{x}\leq ct} u^{\pm}(x,t) = k^{\pm},
$$
and hence
$$
k^- \leq \liminf_{t \to \infty} \inf_{\abs{x}\leq ct} u(x,t) \leq  k^+.
$$
\epf

\section{The characterization of $c^*$ as the slowest speeds of traveling waves}\label{travel}

\subsection{Equivalent integral equations and their upper and lower solutions}\label{upperlower}
In order to establish the existence of travel wave solutions, we fist set up equivalent integral equations.
Similar equivalent integral equations were also used before, see for example,  Wu and Zou \cite{wu2001}, Ma \cite{ma2001,ma2007} and the author \cite{Hwang2009}.
For the convenience of analysis, in this paper and \cite{Hwang2009}, both $\lambda_{1i}, \lambda_{2i}$ are chosen to be positive, and $-\lambda_{1i}, \lambda_{2i}$
are solutions of (\ref{characteis}).

Let $\beta>\max\{\abs{\partial_i f_j(x)}, x \in [0,k^+], i,j=1,...,N\} > 0.$
For $c>c^*$, the two solutions of the following equations,
\begin{equation}\label{characteis}
d_i\lambda^2 -c\lambda -\beta=0, i=1,...,N
\end{equation}
are
$-\lambda_{1i}$ and $\lambda_{2i}$ where
$$
\lambda_{1i}= \frac{-c + \sqrt{c^2+4\beta d_i}}{2d_i}>0, \lambda_{2i}= \frac{c+\sqrt{c^2+4\beta d_i}}{2d_i}>0.
$$
We choose $\beta$ sufficiently large so that
\begin{equation}\label{eq14-1}
\lambda_{2i}>\lambda_{1i} > 2\Lambda_{c}, i=1,...,N.
\end{equation}
Let $u=(u_i) \in \mathcal{C}_{k}$ and define a operator $\mathcal{T}[u]=(\mathcal{T}_i[u])$ by
\begin{equation}\label{eq3}
\begin{split}
\mathcal{T}_i[u](\xi) &= \frac{1}{d_i(\lambda_{1i}+\lambda_{2i})}\Big(\int_{\infty}^{\xi}e^{-\lambda_{1i}(\xi-s)}H_i(u(s))ds\\
&\quad +\int^{\infty}_{\xi}e^{\lambda_{2i}(\xi-s)}H_i(u(s))ds\Big)
\end{split}
\end{equation}
where
$$H_i(u(s))= \beta u_i(s)+ f_i(u(s)),$$
$\mathcal{T}_i[u], i=1,...,N$ is defined on $\mathbb{R}$ if $H_i(u), i=1,2$ is a bounded continuous function. In fact, the following identity
holds
\begin{equation}\label{eq13}
\begin{split}
&\frac{1}{d_i(\lambda_{1i}+\lambda_{2i})}\big (\int_{-\infty}^{\xi}e^{-\lambda_{1i}(\xi-s)} \beta ds +\int^{\infty}_\xi e^{\lambda_{2i}(\xi-s)}\beta ds \big )\\
&=\frac{\beta }{d_i(\lambda_{1i}+\lambda_{2i})} \big (\frac{1}{\lambda_{1i}} + \frac{1}{\lambda_{2i}}\big )= \frac{\beta }{d_i(\lambda_{1i} \lambda_{2i})}\\
&=1.
\end{split}
\end{equation}
We shall show that a fixed point $u$ of $\mathcal{T}$ or solution of the equation
\begin{equation}\label{eq3fixpoint}
u(\xi)= \mathcal{T}[u](\xi)\;\; \xi \in \mathbb{R},
\end{equation}
is a traveling wave solution of (\ref{eq1}) in Lemma \ref{verify}.

\slm\label{verify} Assume $(H1-H2)$ hold.
If $u\in \mathcal{C}_{k}$ is a fixed point of $\mathcal{T}[u]$, $$u(\xi)= \mathcal{T}[u](\xi)\;\; \xi \in \mathbb{R},$$
then
$u$ is a solution of (\ref{eq211}).
\elm
\pf
Note that $H_i(u(s))$ are continuous functions on $\mathbb{R}.$ Thus $\mathcal{T}[u](\xi)$ is defined and differentiable on $\mathbb{R}$.
Direct calculations show
\begin{equation*}
\begin{split}
(\mathcal{T}_i[u](\xi))' &= \frac{1}{d_i(\lambda_{1i}+\lambda_{2i})}\Big(-\lambda_{1i}\int_{-\infty}^{\xi}e^{-\lambda_{1i}(\xi-s)}H_i(u(s))ds\\
&\quad +\lambda_{2i}\int^{\infty}_{\xi}e^{\lambda_{2i}(\xi-s)}H_i(u(s))ds\Big)
\end{split}
\end{equation*}
and
\begin{equation*}
\begin{split}
(\mathcal{T}_i[u](\xi))'' &= \frac{1}{d_i(\lambda_{1i}+\lambda_{2i})}\Big(\lambda_{1i}^2\int_{-\infty}^{\xi}e^{-\lambda_{1i}(\xi-s)}H_i(u(s))ds\\
&\quad +\lambda_{2i}^2\int^{\infty}_{\xi}e^{\lambda_{2i}(\xi-s)}H_i(u(s))ds\\
& \quad -\lambda_{1i}H_i(u(\xi))-\lambda_{2i}H_i(u(\xi)) \Big).
\end{split}
\end{equation*}
Noting that $-\lambda_{1i},\lambda_{2i}$ are solutions of (\ref{characteis}), one can evaluate the following expression
\begin{equation*}
\begin{split}
&(\mathcal{T}_i[u](\xi))''-c(\mathcal{T}_i[u](\xi))' -\beta \mathcal{T}_i[u](\xi) \\
&= \frac{d_i\lambda_{1i}^2+c\lambda_{1i}}{d_i(\lambda_{1i}+\lambda_{2i})}\int_{-\infty}^{\xi}e^{-\lambda_{1i}(\xi-s)}H_i(u(s))ds\\
&\quad + \frac{d_i\lambda_{2i}^2-c\lambda_{2i}}{d_i(\lambda_{1i}+\lambda_{2i})}\int^{\infty}_{\xi}e^{\lambda_{2i}(\xi-s)}H_i(u(s))ds\\
&\quad -H_i(u(\xi))-\beta \mathcal{T}[u](\xi)\\
&=\beta \mathcal{T}_i[u](\xi)-H_i(u(\xi))-\beta \mathcal{T}_i[u](\xi)\\
&=-H_i(u(\xi)).
\end{split}
\end{equation*}
Now if $u(\xi)=\mathcal{T}[u](\xi), \xi \in \mathbb{R}$,  then
$u $ and is a solution of (\ref{eq211}).
\epf

We now define upper and lower solutions of (\ref{eq3fixpoint}), $\phi^+$ and $\phi^-$, which are only continuous on $\mathbb{R}$.
Similar upper and lower solutions have been frequently used in the literatures. See Diekmann \cite{Diekmann1978JMB},
Weinberger \cite{Weinberger1978}, Liu \cite{Lui1989}, Weinberger, Lewis and Li \cite{Weinberger2002-1}, Rass and  Radcliffe \cite{Rass2003},
Weng and Zhao \cite{Weng2006} and more recently, Ma \cite{ma2007}, Fang and Zhao \cite{Fang2009} and Wang \cite{Hwang2009,HwangIntegralDiff}.
In particular, it is believed that the vector-valued lower solutions of the form in this paper  first appeared in \cite{Weng2006}
for multi-type SIS epidemic models. In this paper, the upper and lower solutions here are defined for
general reaction-diffusion systems and we calculate the associated integrals to verify the validity of the upper and lower solutions.

\strdef\label{defupper}
A bounded continuous function $u=(u_i) \in C(\mathbb{R}, [0, \infty)^N)$  is an upper solution of (\ref{eq3fixpoint}) if
$$
\mathcal{T}_i[u](\xi) \leq u_i(\xi), \;\; \text{for all } \xi \in \mathbb{R}, i=1,...,N;
$$
a bounded continuous function $u=(u_i) \in C(\mathbb{R}, [0, \infty)^N)$ is a lower solution of (\ref{eq3fixpoint}) if
$$
\mathcal{T}_i[u](\xi) \geq u_i(\xi), \;\; \text{for all } \xi \in \mathbb{R}, i=1,...,N.
$$
\eeddef

Let $c>c*$ and consider the positive eigenvalue $\Lambda_{c}$ and corresponding eigenvector $\nu_{\Lambda_c}=(\nu^i_{\Lambda_c})$
in Lemma \ref{lmeigen} and $\gamma>1, q>1.$
Define
$$
\phi^+(\xi)=(\phi^+_i),
$$
where
$$
\phi^+_i=\min\{k_i, \nu^i_{\Lambda_c}e^{\Lambda_{c}\xi}\}, i=1,...,N,\; \xi \in \mathbb{R};
$$
and
$$
\phi^-(\xi)=(\phi^-_i),
$$
$$
\phi^-_i=\max\{0, \nu^i_{\Lambda_{c}}e^{\Lambda_{c}\xi}-q\nu^i_{\gamma \Lambda_{c}}e^{\gamma\Lambda_{c}\xi}\},i=1,...,N,\; \xi \in \mathbb{R}.
$$
%

It is clear that if $\xi \geq  \frac{\ln \frac{k_i}{\nu^i_{\Lambda_{c}}}}{\Lambda_{c}}$, $\phi^+_i(\xi)=k_i$, and
$\xi < \frac{\ln \frac{k_i}{\nu^i_{\Lambda_{c}}}}{\Lambda_{c}}$, $\phi^+_i(\xi)=\nu^i_{\Lambda_{c}}e^{\Lambda_{c}\xi},i=1,...,N.$

Similarly, if  $\xi \geq \frac{\ln( q \frac{\nu^i_{\gamma \Lambda_{c}}}{\nu^i_{\Lambda_{c}}}) }{(1-\gamma)\Lambda_{c}}$, $\phi^-_i(\xi)=0$,
and for $\xi < \frac{\ln( q \frac{\nu^i_{\gamma \Lambda_{c}}}{\nu^i_{\Lambda_{c}}}) }{(1-\gamma)\Lambda_{c}}$,
$$\phi^-_i(\xi)=\nu^i_{\Lambda_{c}}e^{\Lambda_{c}\xi}-q\nu^i_{\gamma \Lambda_{c}}e^{\gamma\Lambda_{c}\xi}, i=1,...,N.$$
We choose $q>1$ large enough that
$$
\frac{\ln( q \frac{\nu^i_{\gamma \Lambda_{c}}}{\nu^i_{\Lambda_{c}}}) }{(1-\gamma)\Lambda_{c}}<\frac{\ln \frac{k_i}{\nu^i_{\Lambda_{c}}}}{\Lambda_{c}}, i=1,...,N
$$
and then
$$
\phi^+_i(\xi) > \phi^-_i(\xi), i=1,...,N, \xi \in \mathbb{R}.
$$

\begin{figure}
\begin{center}
  \includegraphics[width=3cm]{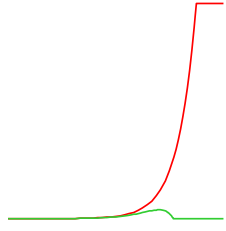}\\
  \caption{For each i, the curve above is $\phi^+_i$ and the below is $\phi^-_i$.}\label{upperlowerslu}
\end{center}
\end{figure}

We now state that $\phi^+$ and $\phi^-$ are upper and lower solution of (\ref{eq3fixpoint}) respectively
and their proofs will be carried out in Appendix through careful analysis of the associated integrals.

\slm\label{upper} Assume $(H1)-(H3)$ hold and (\ref{eq1}) is cooperative.  For any $c > c^*$, $\phi^+$ defined above is an upper solution of (\ref{eq3fixpoint}).
\elm
\slm\label{sub} Assume $(H1)-(H3)$ hold and $\gamma$ satisfies (\ref{gamma}). For any $c > c^*$, $\phi^-$ defined above is a lower solution  of (\ref{eq3fixpoint}) if  $q$ (which is independent of $\xi$) is sufficiently large.
\elm

\subsection{Proof of Theorem \ref{th30}  (iii) when (\ref{eq1}) is cooperative }\label{proofMono}

In this section, we assume that (\ref{eq1}) is cooperative and prove Theorem \ref{th30} (iii). In this case, $f^{\pm}=f$. Many results in this section are standard and
for the verification of continuity and compactness of the operator. See, for example, Ma \cite{ma2001,ma2007} and Wang \cite{Hwang2009}.
Define the following Banach space
$$
\mathcal{E}_{\varrho}= \{u=(u_i): u_i \in C(\mathbb{R}),\sup_{ \xi\in \mathbb{R}} \abs{u_i(\xi)}e^{-\varrho \xi} < \infty, i=1,...,N\}
$$
equipped with weighted norm
$$
\norm{u}_{\varrho}=\sum_{i=1}^N\sup_{ \xi\in \mathbb{R}} \abs{u_i(\xi)}e^{-\varrho \xi},
$$
where $C(\mathbb{R})$ is the set of all continuous functions on $\mathbb{R}$ and $\varrho$ is a positive constant such that $\varrho<\Lambda_{c}.$  It follows that
$\phi^+ \in \mathcal{E}_{\varrho}$ and $\phi^-\in \mathcal{E}_{\varrho}.$
Consider the following set
$$
\mathcal{A}=\{ u=(u_i): u_i \in C(\mathbb{R}) \in \mathcal{E}_{\varrho}, \phi_i^-(\xi) \leq u_i \leq \phi_i^+(\xi), \xi \in \mathbb{R}, i=1,...,N.\}
$$
We shall show the following lemma.

\slm\label{monotone} Assume $(H1)-(H3)$ hold and $\partial_if_j\geq 0, i\neq j$ on $[0,k]$.  Then $\mathcal{T}$ defined in (\ref{eq3}) is monotone
and therefore $\mathcal{T}(\mathcal{A}) \subseteq \mathcal{A}$. Furthermore, $\mathcal{T}_i[u]$ is nondecreasing if $u \in \mathcal{A}$
and all of $u_i$ are nondecreasing.
\elm
\pf  Note that $H_i(u(\xi))$ and $\mathcal{T}[u](\xi)$ are bounded continuous functions on $\mathbb{R}$ if $u \in \mathcal{A}$.
Note  $\beta>\max\{\abs{\partial_i f_i(u)}, u \in [0,k], i=1,...,N\}>0$, $\partial_ig_j(u) \geq 0, u \in [0,k], i \neq j$.
For any $u=(u_i), v=(u_i)\in \mathcal{A}$ with $u_i(\xi) \geq v_i(\xi), \xi \in \mathbb{R},$, we have, for $\xi \in \mathbb{R}$
\begin{equation}\label{eq39}
\begin{split}
&H_i(u(\xi))-H_i(v(\xi))\\
&=\beta(u_i(\xi)-v_i(\xi))+ \int_0^1 \frac{\partial f_i}{\partial u}(su(\xi)+(1-s)v(\xi))ds (u(\xi)-v(\xi))\\
                                   & \geq 0.
\end{split}
\end{equation}

If $u \in \mathcal{A}$ and $u_i$ are nondecreasing, we get, for $i =1,...,N,  \xi\in \mathbb{R}$ and $\xi_1>0$,
\begin{equation}\label{eq40}
\begin{split}
\mathcal{T}_i[u](\xi+\xi_1)-\mathcal{T}_i[u](\xi)&=\frac{1}{d_i(\lambda_{1i}+\lambda_{2i})}\big(\int_{-\infty}^{\xi+\xi_1}e^{-\lambda_{1i}(\xi+\xi_1-s)}H_i(u(s))ds \\
&\quad +\int^{\infty}_{\xi+\xi_1}e^{\lambda_{2i}(\xi+\xi_1-s)}H_i(u(s))ds\\
&\quad -\int_{-\infty}^{\xi}e^{-\lambda_{1i}(\xi-s)}H_i(u(s))ds\\
&\quad -\int^{\infty}_{\xi}e^{\lambda_{2i}(\xi-s)}H_i(u(s))ds\big)\\
&=\frac{1}{d_i(\lambda_{1i}+\lambda_{2i})}\big(\int_{-\infty}^{\xi}e^{-\lambda_{1i}(\xi-s)}H_i(u(s+\xi_1))ds \\
&\quad -\int_{-\infty}^{\xi}e^{-\lambda_{1i}(\xi-s)}H_i(u(s))ds\\
&+\int^{\infty}_{\xi}e^{\lambda_{2i}(\xi-s)}H_i(u(s+\xi_1))ds\\
&-\int^{\infty}_{\xi}e^{\lambda_{2i}(\xi-s)}H_i(u(s))ds\big).
\end{split}
\end{equation}
It follows from (\ref{eq39}) that $\mathcal{T}_i[u](\xi+\xi_1)-\mathcal{T}_i[u](\xi)\geq 0$ for $\xi\in \mathbb{R}$ and $\xi_1>0$.
\epf

Now we shall show that $\mathcal{T}[u]$ is continuous and maps a bounded set in $\mathcal{A}$ into a compact set.

\slm\mlabel{conunity} Assume $(H1)-(H3)$ hold. Then
$\mathcal{T}: \mathcal{A} \rightarrow \mathcal{E}_{\varrho}$ is continuous with the weighted norm $\norm{.}_{\varrho}$.
\elm
\pf Let
$$L=\max\{\abs{\partial_i f_j(u)}, u \in [0,k], i=1,...,N\}.$$
For any $u=(u_i), v=(v_i)\in \mathcal{A}$, we have, for $\xi \in \mathbb{R}$
\begin{equation}\label{eq45}
\begin{split}
&\abs{H_i(u(\xi))-H_i(v(\xi))} e^{-\varrho \xi}\\
& \leq\beta \abs{u_i(\xi)-v_i(\xi)}e^{-\varrho \xi}+|\int_0^1 \frac{\partial f_i}{\partial u}(su(\xi)+(1-s)v(\xi))ds (u(\xi)-v(\xi))|e^{-\varrho \xi}\\
& \leq (\beta +L)\norm{u-v}_{\varrho}\\
\end{split}
\end{equation}
Thus, we obtain
\begin{equation}\label{eq46}
\begin{split}
&\abs{\mathcal{T}_i[u](\xi)-\mathcal{T}_i[v](\xi)} e^{-\varrho \xi}\\
& \leq \frac{1}{(\lambda_{1i}+\lambda_{2i})} \big(\int_{-\infty}^{\xi}e^{-\lambda_{1i}(\xi-s)} \abs{H_i(u(s))-H_i(v(s))}ds\\
                   & \quad +\int_{\xi}^{\infty}e^{\lambda_{2i}(\xi-s)} \abs{H_i(u)(s)-H_i(v)(s)}ds\big)e^{-\varrho \xi}\\
                   &\leq  \frac{(\beta+L) \norm{u-v}_{\varrho}}{(\lambda_{1i}+\lambda_{2i})} \big( \int_{-\infty}^{\xi}e^{-\lambda_{1i}(\xi-s)}e^{\varrho s}ds\\
                   & \quad +\int_{\xi}^{\infty}e^{\lambda_{2i}(\xi-s)}e^{\varrho s}ds\big)e^{-\varrho \xi}\\
                   & = \frac{\lambda_{1i}+\lambda_{2i}}{(\lambda_{1i}+\varrho)(\lambda_{2i}-\varrho)}\frac{(\beta+L)\norm{u-v}_{\varrho}}{(\lambda_{1i}+\lambda_{2i})},
\end{split}
\end{equation}
and
$$
\norm{\mathcal{T}[u]-\mathcal{T}[v]}_{\varrho} \leq \frac{N(\beta+L)}{\min_{i}\{(\lambda_{1i}+\varrho)(\lambda_{2i}-\varrho)\}}\norm{u-v}_{\varrho}.
$$
Thus, $\mathcal{T}[u]$ is continuous.
\epf
\slm\mlabel{compact} Assume $(H1)-(H3)$ hold. Then the set $\mathcal{T}(\mathcal{A})$ is relatively compact in $\mathcal{E}_{\varrho}.$
\elm
\pf
Let $\mathcal{N}_i=\max_{u \in \mathcal{A}, \xi \in \mathbb{R}}H_i(u(\xi))<\infty, i=1,...,N$. Recall that $$\frac{1}{d_i(\lambda_{1i}+\lambda_{2i})}[\int_{-\infty}^te^{-\lambda_{1i}(t-s)}ds +\int^{\infty}_te^{\lambda_{2i}(t-s)}ds]=\frac{1}{\beta}.$$
If $u \in \mathcal{A}$, $\xi \in \mathbb{R}$ and $\delta>0$ ( without loss of generality), we have, $i=1,...,N$
\begin{equation}\label{eq40-1}
\begin{split}
&\mathcal{T}_i[u](\xi+\delta)-\mathcal{T}_i[u](\xi)\\
&=\frac{1}{d_i(\lambda_{1i}+\lambda_{2i})}\big(\int_{-\infty}^{\xi+\delta}e^{-\lambda_{1i}(\xi+\delta-s)}H_i(u(s))ds\\
&\quad +\int^{\infty}_{\xi+\delta}e^{\lambda_{2i}(\xi+\delta-s)}H_i(u(s))ds\\
&\quad -\int_{-\infty}^{\xi}e^{-\lambda_{1i}(\xi-s)}H_i(u(s))ds\\
&\quad -\int^{\infty}_{\xi}e^{\lambda_{2i}(\xi-s)}H_i(u(s))ds\big)\\
&=\frac{1}{d_i(\lambda_{1i}+\lambda_{2i})}\Bigg(\int_{-\infty}^{\xi}e^{-\lambda_{1i}(\xi-s)}\big(e^{-\lambda_{1i}\delta}H_i(u(s))-H_i(u(s))\big)ds\\
&\quad +\int^{\infty}_{\xi}e^{\lambda_{2i}(\xi-s)}\big(e^{\lambda_{2i}\delta}H_i(u(s))-H_i(u(s))\big)ds\\
&\quad + \int_{\xi}^{\xi+\delta}e^{-\lambda_{1i}(\xi+\delta-s)}H_i(u(s))ds\\
&\quad -\int_{\xi}^{\xi+\delta}e^{\lambda_{2i}(\xi+\delta-s)}H_i(u(s))ds\Bigg),
\end{split}
\end{equation}
and
\begin{equation*}
\begin{split}
\mathcal{T}_i[u](\xi+\delta)-\mathcal{T}_i[u](\xi)|& \leq \max\{\abs{e^{-\lambda_{1i} \delta}-1},\abs{e^{\lambda_{2i}\delta}-1}\}\frac{\mathcal{N}_i}{\beta}\\
& \quad + \delta \frac{\mathcal{N}_i}{d_i(\lambda_{1i}+\lambda_{2i})}  + \delta e^{\lambda_{2i} \delta} \frac{\mathcal{N}_i}{d_i(\lambda_{1i}+\lambda_{2i})}.
\end{split}
\end{equation*}
Thus we establish that
\begin{equation}\label{unform}
\lim_{\delta \to 0}(\mathcal{T}_i[u](\xi+\delta)-\mathcal{T}_i[u](\xi))=0, \text{ uniformly for all } u \in \mathcal{A}, \xi\in \mathbb{R}, i=1,...,N.
\end{equation}
Take any sequence $(u^n)=(u_i^n ) \in \mathcal{A}$ and let $v^n=(v_i^n)=\mathcal{T}[u^n]$. From Lemma \ref{monotone} and (\ref{unform}),
$(v^n)$ is uniformly bounded on $\mathbb{R}$ and uniformly equicontinuous.  For $I_m = [-m, m]$, $m \in \mathbb{N}$,
by Ascoli's theorem and the standard diagonal process, we can construct subsequences $(u^{n_m})$ of $(u^n)$ such that there is a function
$v=(v_i), v_i \in C(-\infty, \infty), i=1,...,N$ and $\big(v^{n_m}=\mathcal{T}[u^{n_m}]\big)$ uniformly converges to $v$ on
each $I_m$ for $m \in \mathbb{N}$.  Now we need to show that $v \in \mathcal{A}$ and $\norm{v^{n_m}-v}_{\varrho} \to 0$ as $n_m \to \infty$.
By Lemma \ref{monotone}, $\phi_i^-(\xi) \leq v_i(\xi) \leq \phi_i^+(\xi), i=1,...,N$ for all $\xi  \in \mathbb{R}$, and therefore $v \in \mathcal{A}$.
Note that
$$
\lim_{\xi \to \pm \infty}(\phi_i^+(\xi)-\phi_i^-(\xi))e^{-\varrho \xi}=0, i=1,...,N.
$$
For any $\epsilon>0$, we can find $K_0>0$ such that if $\abs{\xi} > K_0$, then, for all $m \in \mathbb{N}$
$$
\abs{v_i^{n_m}(\xi)-v_i}e^{-\varrho \xi} \leq (\phi_i^+(\xi)-\phi_i^-(\xi))e^{-\varrho \xi} < \epsilon, i=1,...,N.
$$
On the other hand, on $[-I_m, I_m]$, $(v^{n_m})$ uniformly converges to $v$. Thus there exists a $L>0$ such that,
for $n_m>L$
$$
\abs{v_i^{n_m}(\xi)-v_i}e^{-\varrho \xi} < \epsilon,\;\; \xi \in [-K_0, K_0], i=1,...,N.
$$
Consequently, if $n_m>L$, the following inequality is true for all $\xi  \in \mathbb{R}$
$$
\abs{v_i^{n_m}(\xi)-v_i}e^{-\varrho \xi} < \epsilon, i=1,...,N.
$$
Thus $\norm{v^{n_m}-v}_{\varrho} \to 0$ as $n_m \to \infty$.\epf

Now we are in a position to prove Theorem \ref{th30} when (\ref{eq1}) is cooperative.

Define the following iteration
\begin{equation}\label{iteration1}
u^1=(u_i^1)=\mathcal{T}[\phi^+], \;\; u_{n+1}=(u_i^n)=\mathcal{T}[u^n], n>1.
\end{equation}
From Lemmas \ref{upper}, \ref{sub}, \ref{monotone}, $u_n$ is nondecreasing on $\mathbb{R}$ and $$
\phi_i^-(\xi) \leq u_i^{n+1}(\xi) \leq u_i^n(\xi) \leq \phi_i^+(\xi), \xi \in \mathbb{R}, \;n \geq 1,n=1,...,N.
$$
By Lemma \ref{compact} and monotonicity of ($u_n$), there is $u \in \mathcal{A}$ such that $\lim_{n\to \infty}\norm{u_n-u}_{\varrho}=0$.  Lemma \ref{conunity}
implies that $\mathcal{T}[u]=u$. Furthermore, $u$ is nondecreasing. It is clear that $\lim_{\xi \to -\infty}u_i(\xi)=0,i=1,...,N$.
Assume that $\lim_{\xi \to \infty}u_i(\xi)=k_i',i=1,...,N$
$k'_i>0, i=1,...,N$ because of $u \in \mathcal{A}$. Applying the Dominated convergence theorem  to (\ref{eq3}), we get $k_i'=\frac{1}{\beta}(\beta k_i'+f_i(k_1',...,k'_n)$
By (H2), $k_i'=k_i$. Finally, note that
$$
\nu^i_{\Lambda_{c}}(e^{\Lambda_{c} \xi}-q\e^{\gamma \Lambda_{c}\xi}) \leq u_i(\xi) \leq \nu^i_{\Lambda_{c}}e^{\Lambda_{c}\xi}, \xi \in \mathbb{R}.
$$
We immediately obtain
\begin{equation}\label{aystom1}
\lim_{\xi \to -\infty}u_i(\xi)e^{-\Lambda_{c}\xi}=\nu^i_{\Lambda_{c}}, i=1,...,N.
\end{equation}
This completes the proof of Theorem \ref{th30} (iii) when (\ref{eq1}) is cooperative.

\subsection{Proof of Theorem \ref{th30} (iii)}\label{proofofTh1iii}
\pf
Theorem \ref{th30} (iii) is proved when (\ref{eq1}) is cooperative in the last section.
Now we need to prove it in the general case. In order to find traveling waves for (\ref{eq1}), we will apply the Schauder's fixed point theorem.

Let $u=(u_i) \in \mathcal{A}$ and define two integral operators
$$\mathcal{T}^{\pm}[u]=(\mathcal{T}_i^{\pm}[u])$$
 for $f^{-}$ and $f^{+}$
\begin{equation}\label{eq+}
\begin{split}
&\mathcal{T}_i^{\pm}[u](\xi)\\
&=\frac{1}{d_i(\lambda_{1i}+\lambda_{2i})}[\int_{-\infty}^{\xi}e^{-\lambda_{1i}(\xi-s)}H_{i}^{\pm}(u(s))ds +\int^{\infty}_{\xi}e^{\lambda_{2i}(\xi-s)}H_i^{\pm}(u(s))ds]
\end{split}
\end{equation}
and
$$H_i^{\pm}(u(s))= \beta u_i(s)+ f_i^{\pm}(u_(s)).$$
As in Section \ref{proofMono}, both $\mathcal{T}^+$ and $\mathcal{T}^-$ are monotone.  In view of Section \ref{proofMono} and the fact that
$f^-$ is nondecreasing, there exists a nondecreasing fixed point $u^-=(u_i^-)$ of $\mathcal{T}^{-}$ such that
$\mathcal{T}^-[u^-]=u^-$, $\lim_{\xi \to \infty}u_i^{-}(\xi)=k_i^-, i=1,...,N$, and $\lim_{\xi \to -\infty}u_i^{-}(\xi)=0, i=1,...,N$. Furthermore,
$\lim_{\xi \to -\infty}u_i^{-}(\xi)e^{-\Lambda_{c}\xi}=\nu^i_{\Lambda_{c}}, i=1,...,N.$ According to Lemma \ref{upper},
$\phi^+$ (with $k$ being replaced
by $k^{\pm}$) is also a upper solution of $\mathcal{T}^{\pm}$
because the proof of Lemma \ref{upper} is still valid if $f$ is replaced by $f^{\pm}$.  Let
$$
\widetilde{\phi^+}(\xi)=(\widetilde{\phi^+_i}(\xi)),
$$
where
$$
\widetilde{\phi^+_i}(\xi)=\min\{k_i^+, \nu^i_{\Lambda_c}e^{\Lambda_{c}\xi}\}, i=1,...,N,\; \xi \in \mathbb{R};
$$
It follows that $u_i^-(\xi) \leq \widetilde{\phi^+_i}, \xi \in \mathbb{R}, i=1,...,N.$
Now let
\begin{equation}\label{definitionofB}
\mathcal{B}=\{ u: u=(u_i) \in \mathcal{E}_{\varrho}, u_i^-(\xi) \leq u_i(\xi) \leq \widetilde{\phi^+_i}(\xi), \xi \in ( -\infty, \infty), i=1,...,N\},
\end{equation}
where $\mathcal{E}_{\varrho}$ is defined in Section \ref{proofMono}.
It is clear that $\mathcal{B}$ is a bounded nonempty closed convex subset in $\mathcal{E}_{\varrho}$.
Furthermore, we have, for any $u=(u_i) \in \mathcal{B}$
$$
u_i^- = \mathcal{T}_i^-[u^-]\leq \mathcal{T}_i^-[u] \leq \mathcal{T}_i[u] \leq \mathcal{T}_i^+[u]\leq \mathcal{T}_i^+[\widetilde{\phi^+}] \leq \widetilde{\phi^+_i}, i=1,...,N.
$$
Therefore,  $\mathcal{T}: \mathcal{B} \rightarrow \mathcal{B}$. Note that the proof of Lemmas \ref{conunity}, \ref{compact} is valid if (\ref{eq1}) is
not cooperative. In the same way as in Lemmas \ref{conunity}, \ref{compact} , we can show that
$\mathcal{T}: \mathcal{B} \rightarrow \mathcal{B}$ is continuous and maps bounded sets into compact sets.
Therefore, the Schauder Fixed
Point Theorem shows that the operator $\mathcal{T}$ has a fixed point $u$ in $\mathcal{B}$, which is a traveling wave
solution of (\ref{eq1}) for $c>c^*$. Since $u_i^-(\xi) \leq u_i(\xi) \leq \widetilde{\phi^+_i}(\xi), \xi \in ( -\infty, \infty), i=1,...,N$,
it is easy to see that for $i=1,...,N$,
$\lim_{\xi \to -\infty}u_i(\xi)=0$, $\lim_{\xi \to -\infty}u_i(\xi)e^{-\Lambda_{c}\xi}=\nu^i_{\Lambda_{c}}$,
$$k^-\leq \liminf_{\xi \to \infty}u(\xi) \leq \limsup_{\xi \to \infty}u(\xi)\leq k^+$$ and
$0 < u_i^-(\xi) \leq  u_i(\xi) \leq k_i^+, \xi \in ( -\infty, \infty)$.
\epf

\subsection{Proof of Theorem \ref{th30} (iv)} \label{proofofTh1iv}
\pf We adopt the limiting approach in \cite{BrownCarr1977} to prove  Theorem \ref{th30} (iv). For each $n \in \mathbb{N}$,
choose $c_n >c^*$ such that $\lim_{n\to \infty} c_n=c^*.$ According to Theorem \ref{th30} (iii), for each $c_n$
there is a traveling wave solution $u_n=(u^n_i)$ of (\ref{eq1}) such that
$$
u_n=\mathcal{T}[u_n](\xi).
$$
and
$$k^-_i\leq \liminf_{\xi \to \infty}u^n_i(\xi) \leq \limsup_{\xi \to \infty}u^n_i(\xi)\leq k_i^+, i=1,...,N.$$
As it has shown in (\ref{eq40-1}), $(u_n)$ is equicontinuous and uniformly bounded on $\mathbb{R}$, the Ascoli's theorem implies that
there is vector valued continuous function $u=(u^i)$ on $\mathbb{R}$ and subsequence $(u_{n_m})$ of $(u_n)$ such that $$
\lim_{m \to \infty} u_{n_m}(\xi)=u(\xi)
$$
uniformly in $\xi$ on any compact interval of $\mathbb{R}$. Further in view of the dominated convergence theorem we have
$$
u=\mathcal{T}[u](\xi).
$$
Here the underlying $\lambda_{1i}, \lambda_{2i}$ of $\mathcal{T}$ is dependent on $c$ and continuous functions of $c$.
Thus $u$ is a traveling solution of  (\ref{eq1}) for $c=c^*$. Since, for each $c_n$, $u_n \in \mathcal{B}$ where $\mathcal{B}$ is defined
in (\ref{definitionofB}),  it is easy to see that $u$ satisfies
$$k^-\leq \liminf_{\xi \to \infty}u(\xi) \leq \limsup_{\xi \to \infty}u(\xi)\leq k^+$$
Because of the translation invariance of $u_n$, we always can assume that $u_n(0) \leq \frac{1}{2}k^-$ for all $n$.
Consequently $u$ is not a constant traveling solution of (\ref{eq1}).
\epf

\subsection{Proof of Theorem \ref{th30} (v)} \label{proofofTh1v}
\pf
 Suppose, by contradiction, that for some $c \in (0, c^*)$, (\ref{eq1})
has a traveling wave $u(x,t)=u(x + ct)$  with $ \liminf_{\xi \to \infty} u(\xi)>>0$ and $u(-\infty)=0.$  Thus $u(x,t)=u(x + ct)$ can be larger than
a positive vector with arbitrary length. It follows from Theorem \ref{th30} (ii)
$$
\liminf_{t \to \infty} \inf_{\abs{x}  \leq ct} u(x, t) \geq k^->>0, \text{ for } 0< c< c^*
$$
Let $\hat{c} \in (c,c^*)$ and $x=\hat{c}t.$ Then $$
\lim_{t\to \infty} u\big(-(\hat{c}-c)t\big) =\lim_{t\to \infty} u(-\hat{c}t, t) \geq \liminf_{t \to \infty} \inf_{\abs{x}  \leq t\hat{c}} u(x,t) >>0.
$$
However,
$$
\lim_{t \to \infty}u\big(-(\hat{c}-c)t\big)=u(-\infty)=0,
$$
which is a contradiction.
\epf

\section{An example}\label{example}
Weinberger, Kawasaki and Shigesada \cite{WeinbergerKS2009} established the spreading speed for (\ref{eqoper0})
with $h(u_1)$ being a unimodal on $[0,1]$ based on the spreading results for cooperative systems in \cite{Weinberger2002-1}.
Our choice of $h(u_1)$ is slightly different from \cite{WeinbergerKS2009} and simpler than that in \cite{WeinbergerKS2009}.

Our new contribution to (\ref{eqoper0}) is to characterize the spreading speed as the slowest
speed of a family of non-constant traveling wave solutions of (\ref{eqoper0}).
One example of $h(u_1)$  in this paper is $h(u_1)=u_1e^{-u_1}$.
(\ref{eqoper0}) has two equilibriums $(0,0),(0,1)$ and another coexistence equilibrium.  Let $u_1= w_1, u_2=1+w_2$, then (\ref{eqoper0}) can be transformed to
\begin{equation}\label{eqoper00}
\begin{split}
\frac{\partial w_1}{ \partial t}&=d_1 \Delta w_1+w_1[r_1-\alpha -\delta w_1+r_1 w_2]\\
\frac{\partial w_2}{ \partial t}&=d_2\Delta w_2+r_2(1+w_2)[-w_2 +h(w_1)]\\
\end{split}
\end{equation}

In this section, we make the following assumption on $h$.

\begin{enumerate}
  \item[(H4)] \begin{itemize}
               \item[(i)] Assume that $h$ is continuous differentiable on $[0, \infty)$ and $h(0)=0$,
                $h'(0)>0,$ $h(w_1)>0, w_1\in(0,\infty)$. Also assume that $h_m>0$ and $h$ is increasing on $[0,h_m]$, decreases on $[h_m, \infty)$.
                $\lim_{w_1\to \infty} h(w_1)=0.$
               \item[(ii)] Assume that $\frac{h(w_1)}{w_1}$ is strictly decreasing on $(0, \infty)$ and $\lim_{w_1\to \infty}\frac{h(w_1)}{w_1}=0.$
               \item[(iii)]
\begin{equation}\label{eqg77}
\begin{split}
h(w_1)^2+4h(w_1)-4h'(0)w_1 \leq 0, \;\; w_1 \in [0, \infty).
\end{split}
\end{equation}
           \end{itemize}
\end{enumerate}
(H4)(i)(ii) implies that
\begin{equation}\label{eqg779}
\begin{split}
h(w_1) \leq h'(0)w_1, \;\; w_1 \in [0, \infty).
\end{split}
\end{equation}

%
%
We need to verify $h(w_1)=w_1e^{-w_1}$ satisfies (H4). $h(w_1)=w_1e^{-w_1}$ achieves its maximum at $h_m=1$, and is increasing on $[0,h_m]$
and decreasing on $[h_m, \infty)$. In addition, $h'(0)=1$ and $h(w_1)/w_1=e^{-w_1}$ is decreasing for $w_1>0$. It is easy to see that
$e^{x}> x+1, x>0$ and  $e^{-x} < \frac{1}{x+1}, x>0$. Thus, for $w_1>0$
\begin{equation}\label{eq977}
\begin{split}
& h(w_1)^2+4h(w_1)-4h'(0)w_1 \\
&\leq \frac{w_1^2}{2w_1+1}+\frac{4w_1}{w_1+1}-4w_1\\
                     & = \frac{w_1^2(w_1+1)+4w_1(2w_1+1)-4w_1(2w_1+1)(w_1+1)}{(2w_1+1)(w_1+1)}\\
                     & = \frac{-7w_1^3-3w_1^2}{(2w_1+1)(w_1+1)}<0\\
\end{split}
\end{equation}

We also assume that $h_m <k_1$ (otherwise, this problem can be dealt as a cooperative system) and  $\alpha < r_1$. In the nonnegative quadrant, (\ref{eqoper00}) has two equilibrium $(0,0)$ and  $(k_1, k_2)$ satisfying
\begin{equation}\label{equilib2}
\begin{split}
\alpha+\delta k_1  &=r_1+ r_1h(k_1) \\
k_2&=h(k_1).
\end{split}
\end{equation}
We claim that (\ref{equilib2}) has only one positive solution. In fact, the first equation of (\ref{equilib2}) can be rewritten as
\begin{equation}\label{eq8989}
1=\frac{r_1-\alpha+ r_1h(k_1)}{\delta k_1}.
\end{equation}
From (H4)(ii), $\frac{r_1-\alpha+ r_1h(w_1)}{\delta w_1}$ is strictly decreasing
on $(0, \infty)$ and $1=\frac{r_1-\alpha+ r_1h(k_1)}{\delta k_1}$ has only one solution.

In order to use Theorem \ref{th30}, we shall define the two monotone systems. As indicated in Remark \ref{rem1}, similar ideas for constructing
$h^{\pm}$ were used in serval previous works. However the construction of $h^-$ is different.
\begin{equation*}
h^{+}(w_1) = \left\{ \begin{array}{ll}
h(w_1),    & \;\;\;\;\; 0 \leq w_1 \leq h_m, \\[.2cm]
h(h_m), & \;\;\;\;\;  w_1 \geq h_m.
\end{array} \right.
\end{equation*}
and the corresponding cooperative system is
\begin{equation}\label{eqoper1}
\begin{split}
\frac{\partial w_1}{ \partial t}&=d_1 \Delta w_1+w_1[r_1-\alpha -\delta w_1+r_1 w_2]\\
\frac{\partial w_2}{ \partial t}&=d_2\Delta w_2+r_2(1+w_2)[-w_2 +h^+(w_1)]\\
\end{split}
\end{equation}
In a similar manner, one can find (\ref{eqoper1}) has two equilibrium $(0,0)$ and  $(k_1^+, k_2^+)$ satisfying
\begin{equation}\label{equilib2+}
\begin{split}
\alpha+\delta k_1^+  &=r_1+ r_1h^+(k_1^+) \\
k_2^+&=h^+(k_1^+).
\end{split}
\end{equation}
Since  $h^+ \geq  h$, from the first equation of (\ref{eq8989}), it is easily seen that  that $k_1^+ \geq k_1$.
In addition, since $k_1>h_m ,$ we have
$k_2^+=h^+(k_1^+)=h(h_m) \geq h(k_1)=k_2.$

Now there is a $h_0 \in (0, h_m]$ such that $h(h_0)=h(k_1^+)$ and define
\begin{equation*}
h^{-}(w_1) = \left\{ \begin{array}{ll}
h(w_1),    &\;\;\; 0 \leq w_1 \leq h_0, \\[.2cm]
h(k_1^+), & \;\;\; w_1 > h_0.
\end{array} \right.
\end{equation*}
Then $$0<h^-(w_1) \leq h(w_1) \leq  h^+(w_1) \leq h'(0)w_1, w_1 \in (0, k^+_1]$$

\begin{figure}
\begin{center}
  \includegraphics[width=7cm]{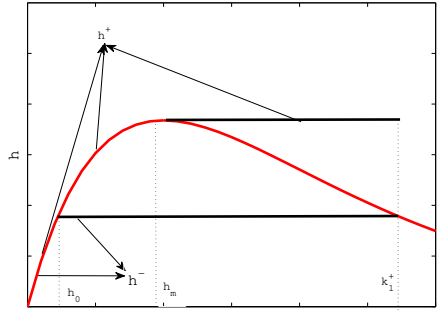}
   \caption{The construction of $h^+$ and $h^-$. The red curve is $h$.}\label{fig1}
\end{center}
\end{figure}

\begin{figure}
\begin{center}
  \includegraphics[width=10cm]{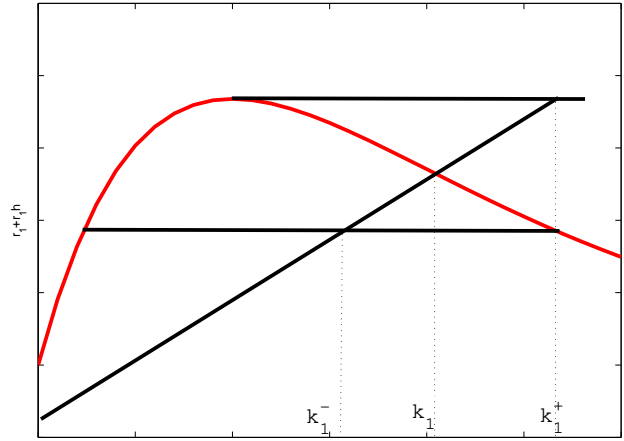}
  \caption{The intersections of $r_1+ r_1h(w_1), r_1+ r_1h^{\pm}(w_1)$ with the line $\alpha+\delta w_1$. The red curve is $r_1+ r_1h$}\label{fig2}
\end{center}
\end{figure}

The corresponding cooperative system for $h^-$ is
\begin{equation}\label{eqoper2}
\begin{split}
\frac{\partial w_1}{ \partial t}&=d_1 \Delta w_1+w_1[r_1-\alpha -\delta w_1+r_1 w_2]\\
\frac{\partial w_2}{ \partial t}&=d_2\Delta w_2+r_2(1+w_2)[-w_2 +h^-(w_1)]\\
\end{split}
\end{equation}
In a similar manner, one can find (\ref{eqoper2}) has two equilibrium $(0,0)$ and  $(k_1^-, k_2^-)$ satisfying
\begin{equation}\label{equilib2-}
\begin{split}
\alpha+\delta k_1^-  &=r_1+ r_1h^-(k_1^-) \\
k_2^-&=h^-(k_1^-).
\end{split}
\end{equation}
Similarly, we have $k_1^- \leq k_1.$  In addition, by the definition of $h^-$, we have
$$k_2^-=h^-(k_1^-) \leq h(k^+_1) \leq  h(k_1)=k_2.$$

Thus, $$
(0,0) <<(k_1^-, k_2^-) \leq (k_1, k_2) \leq (k_1^+, k_2^+).
$$
See Remark \ref{rem3} for (H1)(i). Now it is straightforward to check all other conditions of (H1)(i)-(iv).

The spreading results for the cooperative systems were used to establish
in \cite{WeinbergerKS2009}. We now demonstrate Theorem \ref{th30} can be used to establish spreading speed and
traveling wave solutions of the nonmonotone system (\ref{eqoper00}) and summarize the results in the following theorem.

\sthm\label{th33}  Let $ d_1, \alpha, \delta, r_1, d_2, r_2$ be all positive numbers. Assume $d_1\geq d_2$, $\alpha < r_1$, $k_1>h_m$ and
\begin{equation}\label{eqg_141}
\begin{split}
\delta & \geq  \frac{r_1 r_2 h'(0)}{r_1+r_2-\alpha}.
\end{split}
\end{equation}
Then  the conclusions of Theorem \ref{th30} hold for (\ref{eqoper00}) for the minimum speed $c^*=2\sqrt{(r_1-\alpha)d_1}$,
$\Lambda_c=\frac{c-\sqrt{c^2-4d_1(r_1-\alpha)}}{2d_1}>0$ and  $\nu_{\Lambda_c}$, where $\nu_{\lambda}$ is defined in (\ref{egenvectro4907}).
\ethm
We   now need to check (H2) and (H3). The linearization of (\ref{eqoper00}) at the origin  is
\begin{equation}\label{eqoper10}
\begin{split}
\frac{\partial w_1}{ \partial t}&=d_1 \Delta w_1+(r_1-\alpha)w_1\\
\frac{\partial w_2}{ \partial t}&=d_2\Delta w_2+r_2(h'(0)w_1-w_2)\\
\end{split}
\end{equation}

The matrix in (\ref{egenvalue}) for (\ref{eqoper00}) is
\begin{equation}\label{matrix17}
\begin{split}
A_{\lambda}=(a^{i,j}_{\lambda})=\left(
      \begin{array}{ll}
        d_1\lambda^2+r_1-\alpha & \;\;\;0\\
        r_2h'(0)  & d_2\lambda^2-r_2\\
      \end{array}
    \right)
 \end{split}
\end{equation}
It is easy to see that (H2) holds. In fact, the principle eigenvalue $A_{\lambda}$ is $\Psi(A_{\lambda})=d_1\lambda^2+r_1-\alpha$, which is a convex function of $\lambda$.
And$$\Phi(\lambda)=\frac{\Psi(A_{\lambda})}{\lambda}=\frac{d_1\lambda^2+r_1-\alpha}{\lambda}$$
satisfies the results of Lemma  \ref{lmeigen}. In fact $\Phi(\lambda)$ is  also a strictly convex function of $\lambda$.
The minimum of $\Phi(\lambda)$ is $c^*=2\sqrt{(r_1-\alpha)d_1}$. For each $\lambda>0$, the positive eigenvector of $A_{\lambda}$
corresponding to $\Psi(\lambda)$ is
\begin{equation}\label{egenvectro4907}
\nu_{\lambda}= \left(
                     \begin{array}{c}
                       \nu_{\lambda}^{1} \\
                       \nu_{\lambda}^{2} \\
                     \end{array}
                   \right)
=\left(\begin{array}{c}
                                                                 (d_1-d_2)\lambda^2+r_1+r_2-\alpha\\
                                                                                         r_2h'(0)
                                                                                                        \end{array}
                                                                                                        \right)
\end{equation}
which is also the positive eigenvector of $\frac{1}{\lambda}A_{\lambda}$
corresponding to $\Phi(\lambda)$.

For each $c>c^*$, the left positive solution of $\Phi(\lambda)=c$ is  $$\Lambda_c=\frac{c-\sqrt{c^2-4d_1(r_1-\alpha)}}{2d_1}$$  in Lemma \ref{lmeigen}.
Further from (\ref{egenvectro4907}) we can see that
$$
\frac{\nu_{\lambda}^{2}}{\nu_{\lambda}^{1}}=\frac{r_2h'(0)}{(d_1-d_2)\lambda^2+r_1+r_2-\alpha}=\frac{h'(0)}{\sigma}
$$
where $\sigma=1+\frac{r_1-\alpha +(d_1-d_2)\lambda^2}{r_2}>1.$

In order to verify (H3) for (\ref{eqoper1}) let
$$
(w_1,w_2)=(\theta, \theta \frac{h'(0)}{\sigma})>>(0,0),\; \theta >0.
$$
Thus (H3) is equivalent to the following two inequalities
\begin{equation*}
\begin{split}
w_1[r_1-\alpha -\delta w_1+r_1 w_2]& \leq (r_1-\alpha)w_1\\
r_2(1+w_2)[-w_2 +h^+(w_1)]& \leq r_2(h'(0)w_1-w_2)\\
\end{split}
\end{equation*}
or
\begin{equation}\label{eqoper122}
\begin{split}
\delta w_1& \geq r_1w_2\\
\end{split}
\end{equation}
and
\begin{equation}\label{eqoper123}
\begin{split}
h'(0)w_1+w^2_2& \geq h^+(w_1)(1+w_2)\\
\end{split}
\end{equation}
It follows that the following equality suffices to verify (\ref{eqoper122}):
$$
\delta  \theta  \geq r_1 \theta \frac{h'(0)}{1+\frac{r_1-\alpha}{r_2}},
$$
which is true if (\ref{eqg_141}) holds.

In order to verify (\ref{eqoper123}), following Weinberger, Kawasaki and Shigesada \cite{WeinbergerKS2009}, there always is a positive constant $\varsigma$ such that
$$
w_2=\frac{h'(0)}{\varsigma}w_1.
$$
Substituting $w_2$ into (\ref{eqoper123})  and multiplying $\frac{\varsigma^2}{h'(0)w_1}$ at the both sides, we get
$$
\varsigma^2+h'(0)w_1 \geq \frac{\varsigma^2 h^+(w_1)}{h'(0)w_1}+\varsigma h^+(w_1).
$$
Rearranging the terms produces
\begin{equation}\label{eqoper1233}
\begin{split}
-\varsigma^2\big(1-\frac{h^+(w_1)}{h'(0)w_1}\big)+\varsigma h^+(w_1)-h'(0)w_1 \leq  0, \;\; w_1>0.
\end{split}
\end{equation}
Recall the definition that $h^+(w_1)=h(w_1)$ for $w_1 \leq h_m$ and $h^+(w_1)=h(h_m)$ for $w_1 > h_m$. Since the left side of (\ref{eqoper1233}) is decreasing
in $w_1$ for  $w_1 >h_m$, we only need to verify (\ref{eqoper1233}) for  $w_1 \leq h_m$.

By the quadratic formula, the discriminant of (\ref{eqoper1233}) for $w_1 \leq h_m$ is
$$
h(w_1)^2-4\big(1-\frac{h(w_1)}{h'(0)w_1}\big)h'(0)w_1=h(w_1)^2+4h(w_1)-4h'(0)w_1,
$$
which is negative from (\ref{eqg77}). Thus the left  side of (\ref{eqoper1233}) has no real zeros and
it has to be nonpositive because of the assumption that $1-\frac{h^+(w_1)}{h'(0)w_1}>0$ for $w_1>0$.
In fact, (\ref{eqg77}) is  is one of the three possible conditions in \cite{WeinbergerKS2009}
to guarantee (\ref{eqoper1233}) holds.

Notice that $h^{\pm}$ and $h$ are identical around the origin.  By the exact same arguments (just replacing $k_i^+$ by $k^-_i$ and $h^{+}$ by $h^{-}$), we can verify that (H3) holds for (\ref{eqoper2}) as well.

\section{Appendix}

In this section, we shall verify provide a direct verification of Lemmas \ref{upper} and \ref{sub}. Lower and upper solutions of the equivalent integral
equations (\ref{eq3fixpoint}) play a central role in the construction of fixed points of the equivalent integral equations through monotone iterations.
The lower and upper solutions give asymptotic behavior of traveling wave solutions of (\ref{eq1}).

Wu and Zou \cite{wu2001,wu2008}, Ma \cite{ma2001,ma2007} verify lower and upper solutions through differential equations, and then use them in monotone iterations of equivalent
integral equations. While it was  pointed out in \cite{Boumenir2008} that the upper and lower solutions
for differential equations are required to be smooth for delayed equations, the author \cite{Hwang2009} recently directly verified
$\phi^+$ and $\phi^-$ are indeed lower and upper solutions through the equivalent integral equations for scalar equations,
where the integrals and compare the two sides of (\ref{eq3fixpoint}) were calculated and compared.
Clearly, in this way, the lower and upper solutions are not required to be smooth.

In this appendix, we shall directly verify that,  for $n$-dimensional systems, $\phi^+$ and $\phi^-$ are the lower and
upper solutions of (\ref{eq3fixpoint}). Thus this Appendix can be viewed as a continuation of
\cite{Hwang2009} for the direct verification of non-smooth upper and lower solutions of the equivalent integral equations
for $n$-dimensional systems.

It should be pointed out that the proof of two lemmas in Ma \cite[Lemmas 2.5, 2.6]{ma2001} can significantly simplify the verification of
lower and upper solutions for the equivalent integral equations although the conclusions of the two lemmas in \cite{ma2001} were about lower and upper solutions
for differential equations (see \cite[Section 6.1]{HwangSystemPDE2} for more details). As a result, we always can verify them in a much
simpler way. Nevertheless,  a direct verification can provide a further evidence that $\phi^+$ and $\phi^-$ are lower and upper solutions. In addition, by carefully analyzing
eigenvalues and corresponding eigenvectors, we identify  some identities between the parameters and reveal interesting relations between the parameters.

The results in this appendix are natural extensions of those in the author \cite{Hwang2009} for scalar cases.  As in \cite{Hwang2009},
for $\lambda>0$ let
$$
(M_i(\lambda))=\beta \nu_{\lambda} -\lambda^2D \nu_{\lambda}+A_{\lambda}\nu_{\lambda}
$$
or
\begin{equation}\label{Lambdadef}
M_i(\lambda)=\beta\nu^i_{\lambda} -\nu^i_{\lambda}d_i\lambda^2+\sum_{j=1}^N\nu^j_{\lambda}a^{ij}_{\lambda}, i=1,...,N,
\end{equation}
where $\nu_{\lambda}=(\nu^i_{\lambda})$ is the positive eigenvector of $\frac{1}{\lambda}A_{\lambda}$ in (\ref{egenvalue})
corresponding to the principle eigenvalue $\Phi(\lambda)$.
For $c>c^*$, recall that  $\Phi(\Lambda_{c})=c$. It follows that
$\frac{1}{\Lambda_{c}}A_{\Lambda_{c}} \nu_{\Lambda_c} =\Phi(\Lambda_{c})\nu_{\Lambda_c} =c \nu_{\Lambda_c}$ and
$$
M_i(\Lambda_{c})=(\beta -d_i\Lambda_{c}^2+c \Lambda_{c})\nu^i_{\Lambda_c}, i=1,...,N.
$$
Because of (\ref{eq14-1}),  $M_i(\Lambda_{c})>0,i=1,...,N.$  Noting $M_i(\lambda)$ is continuous with respect to $\lambda$, we always can
choose a $\gamma$ such that
\begin{equation}\label{gamma}
1< \gamma <2, \;\;M_i(\gamma \Lambda_{c})>0, i=1,...,N.
\end{equation}

In order to simply our proofs, we first prove two identities (Lemmas \ref{identity1}, \ref{identity2}),
which are the extension of the identities for scalar cases in \cite{Hwang2009}. Their proofs are almost identical
to those in \cite{Hwang2009} except that eigenvector $\nu^i_{\lambda}$ has to be included.
\slm\label{identity1} Assume $(H1-H2)$ hold. Then for each $c>c^*$
\begin{equation}\label{eq19}
\frac{M_i(\Lambda_{c})}{d_i(\lambda_{1i}+\lambda_{2i})}\Big( \frac{1}{\lambda_{1i}+\Lambda_{c}} +\frac{1}{\lambda_{2i}-\Lambda_{c}} \Big)=\nu^i_{\Lambda_c}, i=1,...,N.
\end{equation}
\elm
\pf
Recall $\lambda_{2i}>\lambda_{1i}> 2\Lambda_{c}>\Lambda_{c}$. It follows  that, $i=1,...,N$,
\begin{equation}\label{eq20}
\begin{split}
&\frac{M_i(\Lambda_{c})}{d_i(\lambda_{1i}+\lambda_{2i})}\Big( \frac{1}{\lambda_{1i}+\Lambda_{c}} +\frac{1}{\lambda_{2i}-\Lambda_{c}} \Big)\\
&=\frac{M_i(\Lambda_{c})}{d_i(\lambda_{1i}+\lambda_{2i})} \frac{(\lambda_{1i}+\lambda_{2i})}{\lambda_{1i}\lambda_{2i} + (\lambda_{2i}-\lambda_{1i})\Lambda_{c} -\Lambda_{c}^2}\\
&=\frac{M_i(\Lambda_{c})}{d_i} \frac{1}{\frac{\beta}{d_i} +\frac{c}{d_i}\Lambda_{c} -\Lambda_{c}^2 }\\
&=M_i(\Lambda_{c}) \frac{1}{\beta +c\Lambda_{c} -d_i\Lambda_{c}^2 }\\
&=\frac{\beta \nu^i_{\Lambda_c} +\nu^i_{\Lambda_c}(c\Lambda_{c} -d_i\Lambda_{c}^2)}{\beta +c\Lambda_{c} -d_i\Lambda_{c}^2 }\\
&=\nu^i_{\Lambda_c}.
\end{split}
\end{equation}

\epf

\slm\label{identity2} Assume $(H1-H2)$ hold and $\gamma$ satisfies (\ref{gamma}). Then for each $c>c^*$, $i=1,...,N$
\begin{equation}\label{eq19-10}
\frac{M_i(\Lambda_{c})}{(\lambda_{1i}+\Lambda_{c})\nu^i_{\Lambda_c}}+\frac{M_i(\Lambda_{c})}{(\lambda_{2i}-\Lambda_{c})\nu^i_{\Lambda_c}}-\frac{M_i(\gamma \Lambda_{c})}{(\lambda_{1i}+\gamma \Lambda_{c})\nu^i_{\gamma\Lambda_c}}-\frac{M_i(\gamma \Lambda_{c})}{(\lambda_{2i}-\gamma \Lambda_{c})\nu^i_{\gamma \Lambda_c}}>0.
\end{equation}
\elm
\pf
Since $\lambda_{2i}>\lambda_{1i}> 2\Lambda_{c}$, it follows that $\lambda_{2i}> \gamma \Lambda_{c}>\Lambda_{c}.$
Lemma \ref{lmeigen} ($\Phi(\gamma\Lambda_{c}) <c$)
implies that, for $i=1,...,N$
\begin{equation}\label{eq360}
\begin{split}
M_i(\gamma\Lambda_{c})&= \beta\nu^i_{\gamma \Lambda_c} -\nu^i_{\gamma\Lambda_c} d_i(\gamma \Lambda_{c})^2+\sum_{j=1,...,N}\nu^j_{\gamma\Lambda_c}a^{ij}_{\gamma\Lambda_{c}}\\
& = \big(\beta -d_i(\gamma \Lambda_{c})^2+\Phi(\gamma\Lambda_{c})\gamma\Lambda_{c}\big)\nu^i_{\gamma\Lambda_c}\\
& < (\beta - d_i(\gamma \Lambda_{c})^2+c\gamma\Lambda_{c})\nu^i_{\gamma\Lambda_c}.
\end{split}
\end{equation}
We also note $M_i(\Lambda_{c})=(\beta+c\Lambda_{c}-d_i\Lambda_{c}^2)\nu^i_{\Lambda_c}$.
Thus, for $i=1,...,N$, we have
\begin{equation}\label{eq36}
\begin{split}
&\frac{M_i(\Lambda_{c})}{(\lambda_{1i}+\Lambda_{c})\nu^i_{\Lambda_c}}+\frac{M_i(\Lambda_{c})}{(\lambda_{2i}-\Lambda_{c})\nu^i_{\Lambda_c}}-\frac{M_i(\gamma \Lambda_{c})}{(\lambda_{1i}+\gamma \Lambda_{c})\nu^i_{\gamma\Lambda_c}}-\frac{M_i(\gamma \Lambda_{c})}{(\lambda_{2i}-\gamma \Lambda_{c})\nu^i_{\gamma \Lambda_c}}\\
&=\frac{(\lambda_{1i}+\lambda_{2i})M_i(\Lambda_{c})}{(\lambda_{1i}\lambda_{2i}+(\lambda_{2i}-\lambda_{1i})\Lambda_{c}-\Lambda_{c}^2)\nu^i_{\Lambda_c}}-\frac{(\lambda_{1i}+\lambda_{2i})M_i(\gamma \Lambda_{c})}{(\lambda_{1i}\lambda_{2i}+(\lambda_{2i}-\lambda_{1i})\gamma \Lambda_{c}-(\gamma \Lambda_{c})^2)\nu^i_{\gamma \Lambda_c}}\\
&=\frac{\frac{\sqrt{c^2+4\beta d_i}}{d_i}M_i(\Lambda_{c})}{(\frac{\beta}{d_i}+\frac{c}{d_i}\Lambda_{c}-\Lambda_{c}^2)\nu^i_{\Lambda_c}}-\frac{\frac{\sqrt{c^2+4\beta d_i}}{d_i}M_i(\gamma \Lambda_{c})}{(\frac{\beta}{d_i}+\frac{c}{d_i}\gamma \Lambda_{c}-(\gamma \Lambda_{c})^2)\nu^i_{\gamma\Lambda_c}}\\
&=\frac{\sqrt{c^2+4\beta d_i}M_i(\Lambda_{c})}{(\beta+c\Lambda_{c}-d_i\Lambda_{c}^2)\nu^i_{\Lambda_c}}-\frac{\sqrt{c^2+4\beta d_i}M_i(\gamma \Lambda_{c})}{(\beta+c\gamma \Lambda_{c}-d_i(\gamma \Lambda_{c})^2)\nu^i_{\gamma\Lambda_c}}\\
&=\sqrt{c^2+4\beta d_i}\big( \frac{M_i(\Lambda_{c})}{(\beta+c\Lambda_{c}-d_i\Lambda_{c}^2)\nu^i_{\Lambda_c}}-\frac{M_i(\gamma \Lambda_{c})}{(\beta +c\gamma \Lambda_{c}-d_i(\gamma \Lambda_{c})^2)\nu^i_{\gamma\Lambda_c}}\big )\\
&=\sqrt{c^2+4\beta d_i}\big(1 -\frac{M_i(\gamma \Lambda_{c})}{(\beta +c\gamma \Lambda_{c}-d_i(\gamma \Lambda_{c})^2)\nu^i_{\gamma \Lambda_c}}\big )\\
& >0.
\end{split}
\end{equation}
This completes the proof.
\epf

\subsection{Proof of Lemma \ref{upper}}
The proof of Lemma \ref{upper} is almost identical to that in \cite{Hwang2009} for the scalar case except that the eigenvector $\nu^i_{\lambda}$ needs to
be included, and delay terms are not present here.
\pf
 Let $\xi^*_i= \frac{\ln \frac{k_i}{\nu^i_{\Lambda_{c}}}}{\Lambda_{c}}, i=1,...,N.$ Then $\phi^+_i(\xi)=k_i$ if $\xi \geq \xi^*_i$, and $\phi^+_i(\xi)=\nu^i_{\Lambda_{c}}e^{\Lambda_{c}\xi}$ if $\xi < \xi^*_i, i=1,...,N.$
Note that $\phi^+_i(\xi) \leq \nu^i_{\Lambda_{c}}e^{\Lambda_{c}\xi}, \xi \in \mathbb{R}.$
In view of (H1)-(H3) we have, for $\xi \leq \xi^*_i$
\begin{equation*}
\begin{split}
H_i(\phi^+(\xi)) &=\beta \phi_i^+(\xi)+ f_i(\phi^+(\xi))\\
& \leq\beta \phi_i^+(\xi)+ \sum_{j=1}^n\partial_j f_i(0)\phi_j^+(\xi)\\
& \leq M_i(\Lambda_{c}) e^{\Lambda_{c}\xi},
\end{split}
\end{equation*}
where $M_i(.)$ is defined in (\ref{Lambdadef}).
For $\xi \geq \xi^*_i$, because of the assumption that (\ref{eq1}) is cooperative and $\beta$ is sufficiently large, we have
\begin{equation*}
\begin{split}
H_i(\phi^+(\xi)) & \leq \beta k_i + f_i(k)\\
               & =\beta k_i.
\end{split}
\end{equation*}
Thus, for $\xi \geq \xi^*_i, i=1,...,N$, we obtain
\begin{equation}\label{eq12}
\begin{split}
\mathcal{T}_i[\phi^+](\xi)& \leq \frac{M_i(\Lambda_{c})}{d_i(\lambda_{1i}+\lambda_{2i})}\int_{-\infty}^{\xi^{*}_i}e^{-\lambda_{1i}(\xi-s)} e^{\Lambda_{c}s}ds\\
                   & \quad +\frac{1}{d_i(\lambda_{1i}+\lambda_{2i})}[\int_{\xi^*_i}^\xi e^{-\lambda_{1i}(\xi-s)} \beta k_ids +\int^{\infty}_\xi e^{\lambda_{2i}(\xi-s)}\beta k_ids].\\
\end{split}
\end{equation}
Thus in view of (\ref{eq13}), we add and subtract the term $\frac{\beta k_i}{d_i(\lambda_{1i}+\lambda_{2i})}\int_{-\infty}^{\xi^{*}_i}e^{-\lambda_{1i}(\xi-s)} ds$ at the left of (\ref{eq12}).
Now for $i=1,...,N, \xi \geq \xi^*_i$, noting that $\nu^i_{\Lambda_{c}}e^{\Lambda_{c}\xi^*_i}=k_i$, (\ref{eq12}) can be written as
\begin{equation}\label{eq14}
\begin{split}
\mathcal{T}_i[\phi^+](\xi)& \leq k_i+\frac{1}{d_i(\lambda_{1i}+\lambda_{2i})} \big (M_i(\Lambda_{c})\int_{-\infty}^{\xi^{*}_i}e^{-\lambda_{1i}(\xi-s)} e^{\Lambda_{c}s}ds\\
&\quad -\beta k_i\int_{-\infty}^{\xi^{*}_i}e^{-\lambda_{1i}(\xi-s)} ds\big )\\
                      & = k_i+\frac{1}{d_i(\lambda_{1i}+\lambda_{2i})} \big (M_i(\Lambda_{c}) \frac{e^{-\lambda_{1i} \xi}e^{(\lambda_{1i}+\Lambda_{c})\xi^*_i}}{\lambda_{1i}+\Lambda_{c}}-\beta k_i\frac{e^{-\lambda_{1i} \xi}e^{\lambda_{1i}\xi^*_i}}{\lambda_{1i}}\big )\\
                      & = k_i+\frac{k_ie^{-\lambda_{1i} \xi}e^{\lambda_{1i} \xi^*_i}}{d_i(\lambda_{1i}+\lambda_{2i})} \big (\frac{M_i(\Lambda_{c})}{(\lambda_{1i}+\Lambda_{c})\nu^i_{\Lambda_{c}}}-\frac{\beta}{\lambda_{1i}}\big )\\
                      & = k_i+\frac{k_ie^{-\lambda_{1i} \xi}e^{\lambda_{1i} \xi^*_i}}{d_i(\lambda_{1i}+\lambda_{2i})(\lambda_{1i}+\Lambda_{c})\lambda_{1i}\nu^i_{\Lambda_{c}}} \big (\lambda_{1i}(M_i(\Lambda_{c})-\nu^i_{\Lambda_{c}}\beta)-\nu^i_{\Lambda_{c}}\beta \Lambda_{c}\big ).
\end{split}
\end{equation}
Since $M_i(\Lambda_{c})=\beta\nu^i_{\Lambda_{c}} -\nu^i_{\Lambda_{c}} d_i\Lambda_{c}^2+c \Lambda_{c}\nu^i_{\Lambda_{c}}, i=1,...,N$, we have, $i=1,...,N$
\begin{equation}\label{eq15}
\begin{split}
&\lambda_{1i}(M_i(\Lambda_{c})-\nu^i_{\Lambda_{c}}\beta)-\beta \nu^i_{\Lambda_{c}}\Lambda_{c}\\
&= \nu^i_{\Lambda_{c}}\lambda_{1i}(c \Lambda_{c}- d_i\Lambda_{c}^2)-\beta \nu^i_{\Lambda_{c}}\Lambda_{c}\\
&= \nu^i_{\Lambda_{c}}\frac{4\beta d_i}{2d_i(\sqrt{c^2+4\beta d_i} + c)}(c \Lambda_{c}- d_i\Lambda_{c}^2)-\beta \nu^i_{\Lambda_{c}}\Lambda_{c}\\
                      & = \nu^i_{\Lambda_{c}}\frac{4\beta (c \Lambda_{c}- d_i\Lambda_{c}^2)-2 (\sqrt{c^2+4\beta d_i} + c) \beta \Lambda_{c}} {2(\sqrt{c^2+4\beta d_i}+c)}\\
                      & = \nu^i_{\Lambda_{c}} \frac{2 \beta\Lambda_{c} \big( 2 c-2d_i\Lambda_{c}-\sqrt{c^2+4\beta d_i} -c \big )} {2(\sqrt{c^2+4\beta d_i}+c)}\\
                      & = \nu^i_{\Lambda_{c}} \frac{2 \beta\Lambda_{c} \big( c-2d_i\Lambda_{c}-\sqrt{c^2+4\beta d_i} \big )} {2(\sqrt{c^2+4\beta d_i}+c)}\\
                      & <0.\\
\end{split}
\end{equation}
Combining (\ref{eq14}) and (\ref{eq15}), we see that
for $\xi \geq \xi^*_i, i=1,...,N$,
\begin{equation}\label{eq16}
\mathcal{T}_i[\phi^+](\xi) \leq k_i.
\end{equation}
Similarly, noting $\nu^i_{\Lambda_{c}}e^{\Lambda_{c} \xi^*_i}=k_i$, one can see that, for $\xi \leq \xi^*_i, i=1,...,N$,
\begin{equation}\label{eq17}
\begin{split}
\mathcal{T}_i[\phi^+](\xi)& \leq \frac{M_i(\Lambda_{c})}{d_i(\lambda_{1i}+\lambda_{2i})} \big(\int_{-\infty}^{\xi}e^{-\lambda_{1i}(\xi-s)} e^{\Lambda_{c}s}ds\\
                   & \quad +\int_{\xi}^{\xi^*_i}e^{\lambda_{2i}(\xi-s)} e^{\Lambda_{c}s}ds\big) +\frac{1}{d_i(\lambda_{1i}+\lambda_{2i})}\int^{\infty}_{\xi^*_i}e^{\lambda_{2i}(\xi-s)}\beta k_is\\
                   &=\frac{M_i(\Lambda_{c})}{d_i(\lambda_{1i}+\lambda_{2i})} \big(\frac{e^{\Lambda_{c}\xi}}{\lambda_{1i}+\Lambda_{c}}+\frac{e^{\Lambda_{c}\xi}}{\lambda_{2i}-\Lambda_{c}}\\
                   & \quad - \frac{e^{\lambda_{2i} \xi}e^{-(\lambda_{2i}-\Lambda_{c})\xi^*_i}}{\lambda_{2i}-\Lambda_{c}}\big) +\frac{\beta k_i}{d_i(\lambda_{1i}+\lambda_{2i})} \frac{e^{\lambda_{2i}\xi}e^{-\lambda_{2i}\xi^*_i}}{\lambda_{2i}}\\
                   & = \frac{e^{\Lambda_{c}\xi}M_i(\Lambda_{c})}{d_i(\lambda_{1i}+\lambda_{2i})}\Big( \frac{1}{\lambda_{1i}+\Lambda_{c}} +\frac{1}{\lambda_{2i}-\Lambda_{c}} \Big)\\
                   & + \frac{M_i(\Lambda_{c}) e^{\lambda_{2i}\xi-\lambda_{2i}\xi^*_i}}{d_i(\lambda_{1i}+\lambda_{2i})}\Big( \frac{-k_i}{(\lambda_{2i}-\Lambda_{c})\nu^i_{\Lambda_{c}}} +\frac{\beta k_i}{M(\Lambda_{c})\lambda_{2i}} \Big).
\end{split}
\end{equation}
Since $M_i(\Lambda_{c})=\beta\nu^i_{\Lambda_{c}} -\nu^i_{\Lambda_{c}} d_i\Lambda_{c}^2+c \Lambda_{c}\nu^i_{\Lambda_{c}}, i=1,...,N$, it is easy to see that, by choosing $\beta$ sufficiently large if necessary,
\begin{equation}\label{eq18}
\begin{split}
&\frac{-k_i}{(\lambda_{2i}-\Lambda_{c})\nu^i_{\Lambda_{c}}} +\frac{\beta k_i}{M_i(\Lambda_{c})\lambda_{2i}}= k_i\frac{(-M_i(\Lambda_{c})+ \nu^i_{\Lambda_{c}}\beta)\lambda_{2i} -\Lambda_{c}\nu^i_{\Lambda_{c}}\beta}{\nu^i_{\Lambda_{c}}(\lambda_{2i}-\Lambda_{c})\lambda_{2i}M_i(\Lambda_{c})}\\
&=k_i\frac{\lambda_{2i}(\nu^i_{\Lambda_{c}} \Lambda_{c}^2 d_i- \nu^i_{\Lambda_{c}} \Lambda_{c}c)-\Lambda_{c}\nu^i_{\Lambda_{c}}\beta}{\nu^i_{\Lambda_{c}}(\lambda_{2i}-\Lambda_{c})\lambda_{2i}M_i(\Lambda_{c})}\\
&=k_i\nu^i_{\Lambda_{c}}\Lambda_{c}\frac{\lambda_{2i}( \Lambda_{c} d_i-  c)-\beta}{\nu^i_{\Lambda_{c}}(\lambda_{2i}-\Lambda_{c})\lambda_{2i}M_i(\Lambda_{c})}\\
&=k_i\nu^i_{\Lambda_{c}}\Lambda_{c}\frac{\frac{c+\sqrt{c^2+4\beta d_i}}{2d_i}( \Lambda_{c} d_i-  c)-\beta}{\nu^i_{\Lambda_{c}}(\lambda_{2i}-\Lambda_{c})\lambda_{2i}M_i(\Lambda_{c})}\\
& \leq 0.
\end{split}
\end{equation}
Combining (\ref{eq19}), (\ref{eq17}) and (\ref{eq18}) leads to,
for $\xi \leq \xi^*_i, i=1,...,N$,
$$
\mathcal{T}_i[\phi^+](\xi) \leq \nu^i_{\Lambda_{c}}e^{\Lambda_{c}\xi}.
$$
And therefore, for $\xi \in \mathbb{R} $,
\begin{equation}\label{eq21}
\mathcal{T}_i[\phi^+](\xi) \leq \phi_i^+(\xi), i=1,...,N.
\end{equation}
This completes the proof of Lemma \ref{upper}.
\epf

We now need the following estimate on $f$, which is an application of the Taylor's Theorem for multi-variable functions. Also see \cite{Hwang2009}.

\slm\label{estimageg} Assume $(H1-H2)$ hold. There exist positive constants $b_{ij}, i,j=1,...,N$ such that
$$
f_i(u)\geq \sum_{j=1}^N\partial_j f_i(0)u_j- \sum_{j=1}^Nb_{ij}(u_j)^2, \;\; u=(u_i), u_i \in [0,k_i], i=1,...,N.
$$
\elm

\subsection{Proof of Lemma \ref{sub}}
Again the proof of Lemma \ref{sub} is almost identical to that in \cite{Hwang2009} for the scalar case except that the eigenvector $\nu^i_{\lambda}$ needs to
be included  and delay terms are not present here.

\pf Let $\xi^*_i=\frac{\ln( q \frac{\nu^i_{\gamma \Lambda_{c}}}{\nu^i_{\Lambda_{c}}}) }{(1-\gamma)\Lambda_{c}}, i=1,...,N.$
If $\xi \geq \xi^*_i $, $\phi^-_i(\xi)=0$,
and for $\xi < \xi^*_i$,
$$\phi^-_i(\xi)=\nu^i_{\Lambda_{c}}e^{\Lambda_{c}\xi}-q\nu^i_{\gamma \Lambda_{c}}e^{\gamma\Lambda_{c}\xi}, i=1,...,N.$$
For $\xi \in \mathbb{R}$, it follows that
\begin{equation*}
\begin{split}
H_i(\phi^-(\xi)) &= \beta \phi_i^-(\xi) + f_i(\phi^-(\xi))\\
               & \geq 0,
\end{split}
\end{equation*}
Thus, for $\xi \geq \xi^*_i$, $$\mathcal{T}_i[\phi^-](\xi) \geq \phi_i^-(\xi), i=1,...,N.$$
We now consider the case $\xi < \xi^*_i$.
It is easy to see that
\begin{equation}\label{eq25}
\begin{split}
\nu^i_{\Lambda_{c}}e^{\Lambda_{c}\xi}&\geq \phi_i^-(\xi) \geq \nu^i_{\Lambda_{c}}e^{\Lambda_{c} \xi}-q \nu^i_{\gamma \Lambda_{c}} e^{\gamma \Lambda_{c}\xi}, \;\; \xi \in \mathbb{R}, i=1,...,N.
\end{split}
\end{equation}
In view of  Lemma \ref{estimageg}, (\ref{eq25}), we have, for $\xi \in \mathbb{R}, i=1,...,N$,
\begin{equation}\label{eq31}
\begin{split}
&H_i(\phi^-(\xi))= \beta \phi_1^-(\xi)+ f_i(\phi^-(\xi))\\
             & \geq \beta \phi_i^-(\xi)+ \sum_{j=1}^n\partial_jf_i(0)\phi_j^-(\xi) -\sum_{j=1}^nb_{ij}(\phi_j^-(\xi))^2\\
             & \geq M_i(\Lambda_{c})e^{\Lambda_{c}\xi} -q M_i(\gamma \Lambda_{c})e^{\gamma \Lambda_{c}\xi}-\widehat{M}_i e^{2 \Lambda_{c}\xi},
\end{split}
\end{equation}
where  $M_i(\cdot)$ is defined in (\ref{Lambdadef}) and
\begin{equation}\label{eq33}
\begin{split}
\widehat{M}_i&=\sum_{j=1}^nb_{ij}(\nu^j_{\Lambda_{c}})^2>0.
\end{split}
\end{equation}
Because of $H_i(\phi^-(\xi)) \geq 0$, the term $\int^{\infty}_{\xi^*}e^{\lambda_{2i}(\xi-s)}H_i(\phi^-(s)) ds$ of $\mathcal{T}_i[\phi^-]$ can be ignored in (\ref{eq34}).  Now we are able to estimate $\mathcal{T}[\phi^-]$
for  $\xi \leq \xi^*, i=1,...,N$
\begin{equation}\label{eq34}
\begin{split}
\mathcal{T}_i[\phi^-](\xi)& \geq \frac{1}{d_i(\lambda_{1i}+\lambda_{2i})}\Big(\int_{-\infty}^{\xi}e^{-\lambda_{1i}(\xi-s)} M_i(\Lambda_{c})e^{\Lambda_{c}s}ds\\
                      & \quad -q \int_{-\infty}^{\xi}e^{-\lambda_{1i}(\xi-s)} M_i(\gamma \Lambda_{c})e^{\gamma\Lambda_{c}s}ds\\
                      &\quad-\widehat{M}_i \int_{-\infty}^{\xi}e^{-\lambda_{1i}(\xi-s)} e^{2\Lambda_{c}s}ds\\
                      & \quad+\int_{\xi}^{\xi^*_i} e^{\lambda_{2i}(\xi-s)}M_i(\Lambda_{c}) e^{\Lambda_{c}s} ds -q \int_{\xi}^{\xi^{*}_i}e^{\lambda_{2i}(\xi-s)} M_i(\gamma \Lambda_{c})e^{\gamma\Lambda_{c}s}ds\\
                      & \quad-\widehat{M}_i \int_{\xi}^{\xi^*_i}e^{\lambda_{2i}(\xi-s)} e^{2\Lambda_{c}s} ds \Big )\\
                      & = \frac{1}{d_i(\lambda_{1i}+\lambda_{2i})} \Big( \frac{M_i(\Lambda_{c})e^{\Lambda_{c}\xi}}{\lambda_{1i}+\Lambda_{c}}-q \frac{M_i(\gamma \Lambda_{c})e^{\gamma \Lambda_{c}\xi}}{\lambda_{1i}+\gamma \Lambda_{c}}\\
                      & \quad-\widehat{M}_i\frac{e^{2 \Lambda_{c}\xi}}{\lambda_{1i}+2 \Lambda_{c}}+\frac{e^{\Lambda_{c}\xi^*_i-\lambda_{2i}\xi^*_i +\lambda_{2i}\xi}-e^{\Lambda_{c}\xi}}{\Lambda_{c}-\lambda_{2i}}M_i(\Lambda_{c})\\
                      & \quad-q\frac{e^{\gamma \Lambda_{c}\xi^*_i-\lambda_{2i}\xi^*_i +\lambda_{2i}\xi}-e^{\gamma \Lambda_{c}\xi}}{\gamma \Lambda_{c}-\lambda_{2i}}M_i(\gamma \Lambda_{c})\\
                      & \quad -\widehat{M}_i\frac{e^{2 \Lambda_{c}\xi^*_i-\lambda_{2i}\xi^*_i +\lambda_{2i}\xi}-e^{2\Lambda_{c}\xi}}{2 \Lambda_{c}-\lambda_{2i}} \Big)\\
\end{split}
\end{equation}
In view of the identity (\ref{eq19}), we subtract two terms to make up a term $-\nu^i_{\gamma\Lambda_{c}}qe^{\gamma \Lambda_{c}\xi}$ and thus we need to add the terms. Recall that $\gamma \Lambda_{c}< 2\Lambda_{c} < \lambda_{2i}.$ We ignore two positive
terms $q\frac{M_i(\gamma \Lambda_{c})}{(\lambda_{2i}-\gamma \Lambda_{c})d_i(\lambda_{1i}+\lambda_{2i})}e^{(\gamma \Lambda_{c}-\lambda_{2i})\xi^*+\lambda_{2i}\xi}$ and\\
$\frac{\widehat{M}_i}{(\lambda_{2i}-2 \Lambda_{c})d_i(\lambda_{1i}+\lambda_{2i})}e^{(2\Lambda_{c}-\lambda_{2i})\xi^*+\lambda_{2i}\xi}$.  Thus,
\begin{equation}\label{eq35}
\begin{split}
\mathcal{T}_i[\phi^-](\xi)& \geq \frac{M_i(\Lambda_{c})}{d_i(\lambda_{1i}+\lambda_{2i})} \Big (\frac{1}{\lambda_{1i}+\Lambda_{c}} +\frac{1}{\lambda_{2i}-\Lambda_{c}}\Big)e^{\Lambda_{c}\xi}\\
                      & \quad  -\frac{M_i(\Lambda_{c})}{d_i(\lambda_{1i}+\lambda_{2i})} \Big (\frac{1}{\lambda_{1i}+\Lambda_{c}} +\frac{1}{\lambda_{2i}-\Lambda_{c}}\Big)\frac{\nu^i_{\gamma \Lambda_{c}}}{\nu^i_{\Lambda_{c}}}qe^{\gamma \Lambda_{c}\xi}\\
                      & \quad +\frac{e^{\gamma \Lambda_{c}\xi}}{d_i(\lambda_{1i}+\lambda_{2i})}\Bigg\{q\nu^i_{\gamma \Lambda_{c}}\Bigg(\frac{M_i(\Lambda_{c})}{(\lambda_{1i}+\Lambda_{c})\nu^i_{\Lambda_{c}}}+\frac{M_i(\Lambda_{c})}{(\lambda_{2i}-\Lambda_{c})\nu^i_{\Lambda_{c}}}\\
                      & \quad -\frac{M_i(\gamma \Lambda_{c})}{(\lambda_{1i}+\gamma \Lambda_{c})\nu^i_{\gamma \Lambda_{c}}}-\frac{M_i(\gamma \Lambda_{c})}{(\lambda_{2i}-\gamma \Lambda_{c})\nu^i_{\gamma \Lambda_{c}}}\Bigg)\\
                      & \quad -\frac{\widehat{M}_i}{(\lambda_{1i}+2 \Lambda_{c})}e^{(2-\gamma)\Lambda_{c}\xi}-\frac{M_i(\Lambda_{c})e^{( \Lambda_{c}-\lambda_{2i})\xi^*}}{(\lambda_{2i}-\Lambda_{c})}e^{(\lambda_{2i}-\gamma \Lambda_{c})\xi}\\
                      & \quad -\frac{\widehat{M}_i}{(\lambda_{2i}-2 \Lambda_{c})}e^{(2-\gamma)\Lambda_{c}\xi}\Bigg\}.
\end{split}
\end{equation}
For $\xi \leq \xi^*_i,$ $e^{(2-\gamma)\Lambda_{c}\xi},e^{(\lambda_{2i}-\gamma \Lambda_{c})\xi}$ are bounded above. Because of the identity (\ref{eq19}),
(\ref{eq35}) can be further simplified as, $i=1,...,N$
\begin{equation}\label{eq35.5}
\begin{split}
\mathcal{T}_i[\phi^-](\xi)                      & \geq \nu^i_{\Lambda_{c}}e^{\Lambda_{c}\xi}- q\nu^i_{\gamma \Lambda_{c}}e^{\gamma \Lambda_{c}\xi}\\
                      & \quad +\frac{e^{\gamma \Lambda_{c}\xi}}{d_i(\lambda_{1i}+\lambda_{2i})}\Bigg\{q\nu^i_{\gamma \Lambda_{c}}\Bigg(\frac{M_i(\Lambda_{c})}{(\lambda_{1i}+\Lambda_{c})\nu^i_{\Lambda_{c}}}+\frac{M_i(\Lambda_{c})}{(\lambda_{2i}-\Lambda_{c})\nu^i_{\Lambda_{c}}}\\
                      & \quad -\frac{M_i(\gamma \Lambda_{c})}{(\lambda_{1i}+\gamma \Lambda_{c})\nu^i_{\gamma \Lambda_{c}}}-\frac{M_i(\gamma \Lambda_{c})}{(\lambda_{2i}-\gamma \Lambda_{c})\nu^i_{\gamma \Lambda_{c}}}\Bigg)\\
                      & \quad -\frac{\widehat{M}_i}{(\lambda_{1i}+2 \Lambda_{c})}e^{(2-\gamma)\Lambda_{c}\xi^*_i}-\frac{M_i(\Lambda_{c})e^{( \Lambda_{c}-\lambda_{2i})\xi^*_i}}{(\lambda_{2i}-\Lambda_{c})}e^{(\lambda_{2i}-\gamma \Lambda_{c})\xi^*_i}\\
                      & \quad -\frac{\widehat{M}_i}{(\lambda_{2i}-2 \Lambda_{c})}e^{(2-\gamma)\Lambda_{c}\xi^*_i}\Bigg\}.
\end{split}
\end{equation}

Finally, from (\ref{eq35.5}) and Lemma \ref{identity2}, we conclude that there exists $q>0$, which is independent of $\xi$, such that, for $\xi \leq \xi^*_i$
and $i=1,...,N$
\begin{equation}\label{eq37}
\begin{split}
\mathcal{T}_i[\phi^-](\xi)                      & \geq \nu^i_{\Lambda_{c}}e^{\Lambda_{c}\xi}-q\nu^i_{\gamma\Lambda_{c}} e^{\gamma \Lambda_{c}\xi}.
\end{split}
\end{equation}
And therefore, for $i=1,...,N$
\begin{equation*}
\begin{split}
\mathcal{T}_i[\phi^-](\xi)                      & \geq \phi_i^-(\xi), \;\; \xi \in \mathbb{R}.
\end{split}
\end{equation*}
This completes the proof.
\epf

\end{document}